\documentclass[12pt,reqno]{article}

\usepackage[usenames]{color}
\usepackage{amssymb}
\usepackage{graphicx}
\usepackage{amscd}

\usepackage[colorlinks=true,
linkcolor=webgreen,
filecolor=webbrown,
citecolor=webgreen]{hyperref}

\definecolor{webgreen}{rgb}{0,.5,0}
\definecolor{webbrown}{rgb}{.6,0,0}

\usepackage{color}
\usepackage{fullpage}
\usepackage{float}

\usepackage{graphics,amsmath,amssymb}
\usepackage{amsthm}
\usepackage{amsfonts}
\usepackage{latexsym}

\begin{document}

\vspace*{2.1cm}

\theoremstyle{plain}
\newtheorem{theorem}{Theorem}
\newtheorem{corollary}[theorem]{Corollary}
\newtheorem{lemma}[theorem]{Lemma}
\newtheorem{proposition}[theorem]{Proposition}
\newtheorem{obs}[theorem]{Observation}
\newtheorem{claim}[theorem]{Claim}

\theoremstyle{definition}
\newtheorem{definition}[theorem]{Definition}
\newtheorem{example}[theorem]{Example}
\newtheorem{remark}[theorem]{Remark}
\newtheorem{conjecture}[theorem]{Conjecture}
\newtheorem{question}[theorem]{Question}

\begin{center}

\vskip 1cm

{\Large\bf From semi-total to equitable total colorings, and their induced partitions into perfect and total-perfect codes} %
\vskip 5mm
\large
 Italo J. Dejter

University of Puerto Rico

Rio Piedras, PR 00936-8377

\href{mailto:italo.dejter@gmail.com}{\tt italo.dejter@gmail.com}
\end{center}


\begin{abstract} 
Independently posed by Behzad and Vizing, the Total Coloring Conjecture asserts that the total chromatic number $\chi''(G)$ of a connected simple graph $G$ is either $\Delta(G)+1$ or $\Delta(G)+2$, where $\Delta(G)$ is the largest degree of a vertex of $G$.
 Deciding whether a cubic graph $G$ has $\chi''(G)=\Delta(G)+1$, even if $G$ is bipartite, is NP-hard. The problem persists even for equitable total colorings.
Williams and Holroyd gave a condition for total colorings via the introduction of semi-total colorings. In the present work, focus is set mainly on obtaining, by means of a  variation of Kempe's 1879 graph-coloring algorithm, equitable total colorings out of semi-total colorings of: symmetric cubic graphs (3-cube, Petersen, Heawood, Pappus, Desargues, etc., graphs) and cage graphs (Robertson graph, etc.)
Moreover, total colorings of the 3-cube and subsequent $4k$-polygon prisms ($k>1$) are shown to exist that have perfect codes as their vertex color classes, as well as semi-total colorings of M\"obius ladders that have total-perfect codes as their vertex color classes, only $K_{3,3}$ of which equitable. In addition, a semi-total coloring of the 90-vertex Foster graph is shown to have a vertex color class being a total-perfect code. On the other hand, a semi-total coloring of the 42-vertex $(5,6)$-cage yields a decomposition into 3 copies of the Heawood graph, by the edges of same color endvertices, and 7 6-cycles. Among other cases, non-equitable total colorings of cubic graph expansions are taken to equitable semi-total colorings.
\end{abstract}

\section{A  variation of Kempe's coloring algorithm}\label{s0}

Let $\Delta(G)$ be the largest degree of any vertex of a connected simple graph $G$. The total coloring (or {\it TC}) conjecture \cite{B1,B2,V} asserts that the total chromatic number $\chi''(G)$ of $G$ (see Definition~\ref{tcol}) is either $\Delta(G)+1$ or $\Delta(G)+2$. 
To decide whether a cubic graph $G$ has $\chi''(G)=\Delta(G)+1$, even if $G$ is bipartite cubic, is NP-hard \cite{Dantas,Arroyo}. The problem persists even for equitable total colorings, namely those whose color-class cardinalities differ at most by 1; see Subsection~\ref{sub3}.

Williams and Holroyd \cite{WH} gave a condition to solve TC problems departing from {\it semi-total colorings} (or {\it sTC}s); see
Definition~\ref{acanomas}. 
A  variation of Kempe'a 1879 graph-coloring algorithm~\cite{Kempe1} (that has applications in Register Allocation \cite{Kempe2}) is introduced in Theorem~\ref{t1}, taking semi-total colorings of connected simple graphs and possibly modifying them into equitable total colorings, or in their defect, into equitable semi-total colorings. Theorem~\ref{t1} and Theorem~\ref{fold}, this one dealing with covering graph maps, are applied to three families of graphs, namely symmetric cubic graphs (Section~\ref{s2}), cage graphs (Examples~\ref{e1}, \ref{pete}, \ref{K33} and Section~\ref{last}) and vertex expansions of cubic graphs (Section~\ref{s3}). 

It is to be noted that $4k$-polygon prisms with $k\ge 1$ have TCs whose color vertex classes are perfect codes (Theorem~\ref{girth}) and that M\"obius ladders have sTCs whose color vertex classes are total-perfect codes via all edges with same color endvertices (Theorem~\ref{valentino} dealing with a generalization of M\"obius ladders and Corollary~\ref{ladder}), where the smallest such ladder, namely the utility graph $K_{3,3}$, is the only equitable one of them. In addition, the 90-vertex Foster graph has an sTC just one of whose color classes is a total-perfect code (Theorem~\ref{fost}). Moreover,
an sTC of the $(5,6)$-cage yields its decomposition into three copies of the Heawood graph via all edges with same color endvertices and seven 6-cycles with all remaining edges (Theorem~\ref{lasthm}).
We start with some basic definitions. 

\begin{figure}[htp]
\includegraphics[scale=0.56]{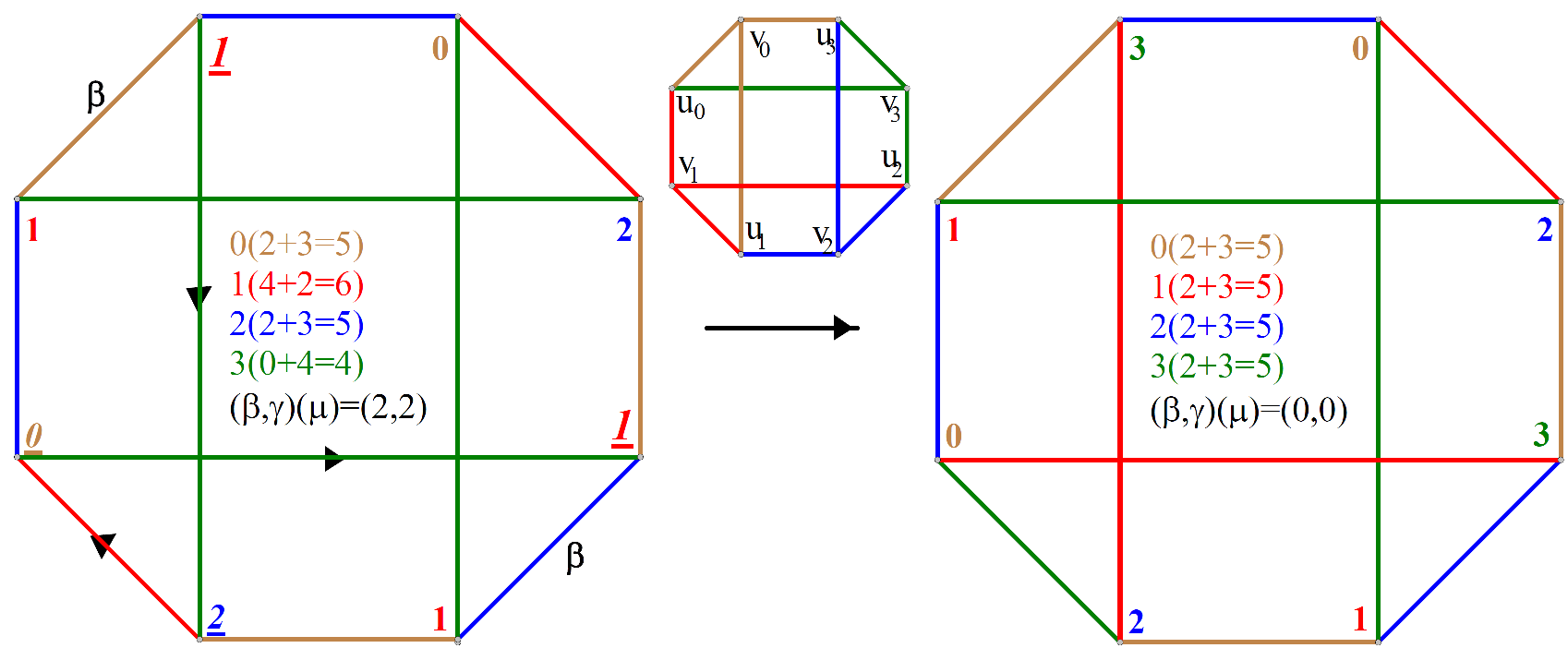}
\caption{From a lacunar sTC to an EqTC for the 3-cube via $\beta$-reduction (Example~\ref{q3}).}
\label{cubo}
\end{figure}

\subsection{Total and semi-total colorings}\label{sub1}

\begin{definition}\label{tcol} Let $G=(V(G),E(G))$ be a connected simple graph. A {\it total coloring} (or {\it TC}) of $G$ is a color assignment to the elements of $G$ (that is: vertices and edges) such that no two incident or adjacent elements are assigned the same color.\end{definition}

\begin{definition}
The {\it total chromatic number} $\chi''(G)$ of a simple connect graph $G$ is the least number of colors needed to color the elements of $G$ such that no two adjacent or incident elements are assigned the same color.
\end{definition}

 A recent survey \cite{tc-as} contains the most up to date resources on TCs. 
 Recent results on the TC Conjecture are found in \cite{Dantas,Mazzu,Feng,Rosen,Vi} including that the total chromatic number of cubic graphs is either 4 or 5 \cite{Feng}. 

\begin{definition} A connected simple graph $G$ is of {\it type} 1 if $\chi''(G)=\Delta(G)+1$ and is of {\it type} 2 if $\chi''(G)>\Delta(G)+1$, \cite{Rosen,V}. 
 \end{definition}

\begin{figure}[htp]
\includegraphics[scale=0.56]{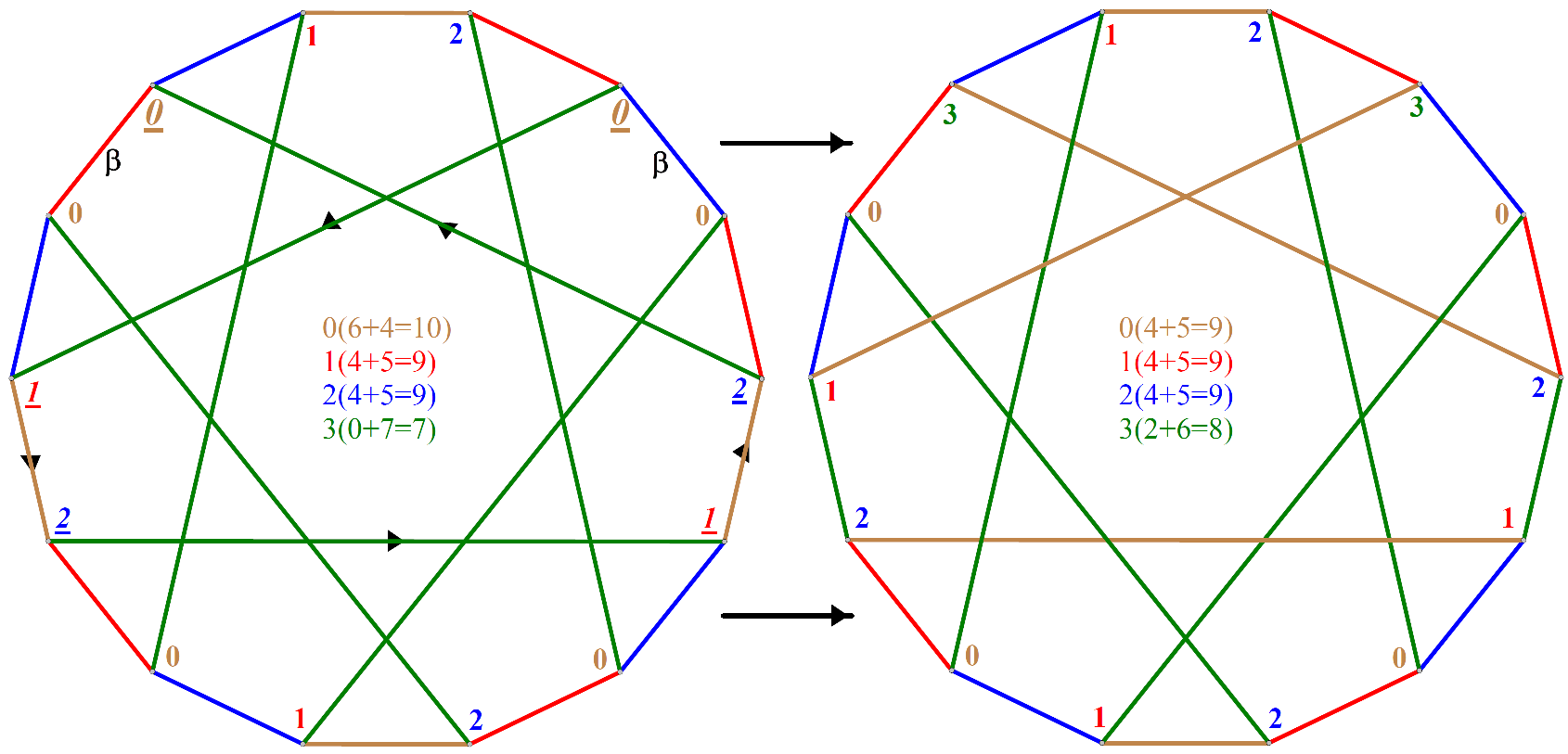}
\caption{From a lacunar sTC to an EqTC for the Heawood graph (Example~\ref{e1}).}
\label{figshp}
\end{figure}

Yap \cite{Yap} showed that: {\bf(a)} a cycle $C_n$ is of type 1 if and only if $n$ is a multiple of 3; {\bf(b)} a complete graph $K_n$ is of type 1 if and only if $n$ is odd.


\begin{definition}\label{acanomas} Given a connected simple graph $G$ with maximum degree $\Delta=\Delta(G)$, a {\it semi-total coloring} (or {\it sTC}) $\mu$ of $G$ is a function $\mu:V(G)\cup E(G)\rightarrow\chi(\Delta)=\{1,2,\ldots,\Delta +1\}$, where $\chi(\Delta)$ is the assumed set of colors for $G$, such that every two adjacent edges of $G$ have distinct colors of $\chi(\Delta)$ and every vertex of $G$ has a color in $\chi(\Delta)$ distinct from the colors of its incident edges.\end{definition} 

The concept of an sTC is similar to that of a TC, except that it is not required that adjacent vertices have distinct colors. The need of using $\Delta+1$ colors to obtain an sTC for $G$ emerges because at least $\Delta+1$ colors are needed, and by Vizing's Theorem~\cite{V}, $G$ has an edge coloring with at most $\Delta+1$ colors, so that there is at least one color free for each vertex.

\begin{definition} Given a connected simple graph $G$ and an sTC $\mu$ of $G$, an edge $v_0v_1$ of $G$ such that $\mu(v_0)=\mu(v_1)$ is said to be a {\it $\beta$-edge} of $G$ {\it with respect to} (or {\it wrt}) $\mu$. The {\it $\beta$-set of $G$ wrt $\mu$} is the set of such $\beta$-edges.
 The number of $\beta$-edges of $G$ wrt $\mu$ is denoted by $\beta(\mu)$.
\end{definition}

\begin{definition} Given a connected simple graph $G$ and an sTC $\mu$ of $G$, the minimum value of $\beta(\mu)$ over all sTCs $\mu$ of $G$ is said to be the {\it beta parameter} $\beta(G)$ of $G$.
\end{definition}

 Williams and Holroyd \cite{WH} showed that: {\bf(a)} if $n\equiv 1$ or 2 (mod 3), then $\beta(C_n)=2$; {\bf(b}) if $n$ is even, then $\beta(K_n)=n/2$.

\subsection{Equitable total coloring}\label{sub3}

\begin{definition} Let $G$ be a connected simple graph and let $\mu$ be an sTC of $G$. The {\it $\mu$-class} of a color $i\in\chi(\Delta)$ is the subset of $V(G)\cup E(G)$ having its elements assigned color $i$ by $\mu$.
\end{definition}

\begin{figure}[htp]
\hspace*{2.0cm}
\includegraphics[scale=0.95]{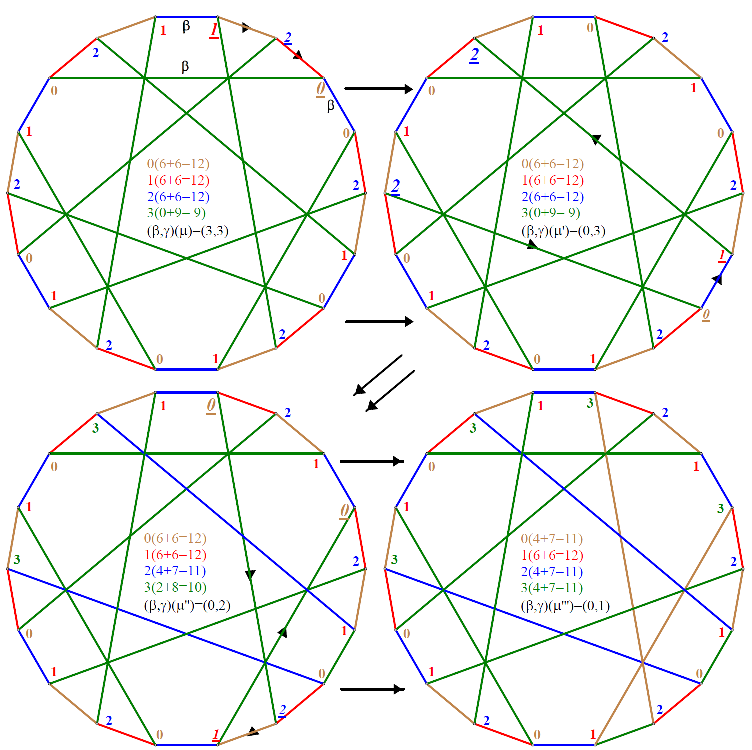}
\caption{From a lacunar TC to an EqTC for the Pappus graph (Example~\ref{Pap}).}
\label{figeqp}
\end{figure}

\begin{definition} 
 An sTC $\mu$ of a connected simple graph $G$ is {\it equitable} if the cardinalities of every two $\mu$-classes differ by at most 1, in which case $\mu$ is {\it equitable sTC}, ({\it EqTC}, if $\mu$ is TC). 
\end{definition}

\begin{definition} Given an sTC $\mu$  of a connected simple graph $G$,
let $\gamma(\mu)$ be the difference between the maximum and minimum cardinalities of $\mu$-classes of $G$.
\end{definition}

\begin{definition} Given a connected simple graph $G$, the minimum value of $\gamma(\mu)$ over all sTCs $\mu$ of $G$ with $\chi(\Delta(G))$ colors is said to be the {\it gamma parameter} $\gamma(G)$ of $G$. 
\end{definition}

It is conjectured that the equitable total chromatic number of a connected simple graph $G$ is at most $\Delta(G) + 2$ (similar to the TC Conjecture) and proved for subcubic graphs \cite{Wang}. Furthermore, it has been shown in \cite{Gui} that, for every subcubic graph $G$, there exists an equitable $k$-total coloring of $G$ for each $k>\Delta +2$, which proves the validity of a conjecture of Fu~\cite{Fu} in the case of subcubic graphs, namely that every $G$ has an equitable total $k$-coloring for each $k\ge max\{\chi''(G),\Delta + 2\}$.

\subsection{Beta and Gamma reductions}\label{sub4}

 With the aim of reducing either the value
$\beta(\mu)>0$ and/or the value $\gamma(\mu)>1$ of an sTC $\mu$ of a connected simple graph $G$, we will employ   
the concept of {\it skeleton} of a connected simple graph $G$ having an sTC $\mu$. 
 We start by presenting Definitions~\ref{dej1}-\ref{ske}.  

\begin{definition}\label{dej1} Let $G$ be a connected simple graph with maximum vertex degree $\Delta=\Delta(G)$. Let $\mu:V(G)\cup E(G)\rightarrow\chi(\Delta)$ be an sTC. Let  $S=(v_0,e_1,v_1,e_2,v_2,\ldots,e_k,v_k)$ be a  path in $G$ on $k+1$ vertices, so its length is $k$ and $S\equiv P_{k+1}$. 
Assume that $\{c_0,c_1\}\in\chi(\Delta)$, where $c_0\neq c_1$. Let $\mu(v_0)=\mu(e_2)=\cdots=\mu(e_{2i})=\cdots=c_0$ and  $\mu(e_1)=\mu(e_3)=\cdots=\mu(e_{2i-1})=\cdots=c_1,$ where  $0<i\le\lceil\frac{k}{2}\rceil$ and $\mu(e_k)\neq\mu(v_k)\in\{c_0,c_1\}$, 
 so that $$\mu_S=(\mu(v_0),\mu(e_1),\mu(e_2),\ldots,\mu(e_{k-1}),\mu(e_k),\mu(v_k))$$ is a 2-valued alternating sequence. 
Then, we say that $S$ is a {\it color alternating path wrt to $\mu$}.
\end{definition} 

\begin{definition} A color alternating path of $G$ wrt $\mu$ is a {\it maximal color alternating path} (or {\it MCAP}) of $G$ wrt $\mu$ if no path in $G$ having $S$ as a shorter subpath is a color alternating path of $G$ wrt $\mu$. 
\end{definition}

\begin{definition}\label{ske} Given an MCAP $S=\{v_0,e_1,v_1,e_2,v_2,\ldots,e_k,v_k\}$ of $G$ wrt $\mu$, the sequence $\rho(S)=(v_0,e_1,e_2,\ldots,e_k,v_k)$ whose terms consists of the two terminal degree-one vertices $v_0$ and $v_k$ and every edge $e_i$ of $S$ in between, for $i=1,\ldots,k$, is said to be the {\it skeleton} of $S$ wrt to $\mu$.
\end{definition}

\begin{theorem}\label{t1} 
 Let $\mu$ be an sTC of a connected simple graph $G$ with maximum degree $\Delta$, where $\beta(\mu)>0$, or\, $\gamma(\mu)>1$. Let $S=\{v_0,e_1,v_1,e_2,v_2,\ldots,e_k,v_k\}$ be an MCAP of $G$ wrt $\mu$, where $\mu(S)=\{c_0,c_1\}$ and $\mu(v_0)=c_0$. Then, there exists an sTC $\mu':V(G)\cup E(G)\rightarrow\chi(\Delta)$ with $\mu'(S)=\{c_0,c_1\}$ and $\mu'$ differing from $\mu$ only on the skeleton $\rho(S)=\{v_0,e_1,\ldots,e_k,v_k\}$, where each value of $\mu'$ differs from the corresponding value of $\mu$. Moreover: {\bf(a)}
 If the set of $\beta$-edges of $\mu$ strictly contains the set of $\beta$-edges of $\mu'$, then $\beta(\mu')<\beta(\mu)$; {\bf(b)} the difference between the maximum and minimum cardinalities is less for the $\mu'$-classes than for the $\mu$-classes if and only if $\gamma(\mu')<\gamma(\mu)$.
\end{theorem}

\begin{example}\label{joder3} Here is a counterexample to the converse of conclusion (a) in Theorem~\ref{t1}.
Take the tree with edges $uv,ua,ub,vc$ and vertex colors $\mu(u)=\mu(v)=\mu(a)=\mu(v)=1$ and $\mu(c)=2$ and assign to those four edges the colors $2,3,4,3$, respectively. Swapping colors 1 and 2  on the maximal single-edge skeleton $\{u,uv,v\}$ preserves the sTC condition but changes the the $\beta$-set from $\{uv,us,ub\}$ to $\{uv,vc\}$, so $\beta$ decreases from 3 to 2 while a new $\beta$-edge appears.
\end{example}

\begin{proof}
First, assume that $k$ is odd. Then we define $\mu'$ by exchanging the values $c_0$ and $c_1$ in $\mu(\rho(S))$ by $c_1$ and $c_0$, respectively, as shown in the following display:

\begin{eqnarray}\label{joder1}\begin{array}{||c||c|c|c|c|c|c|c|c||}\hline\hline
   G&v_0&e_1&e_2&e_3&\ldots&e_{k-1}&e_k&v_k\\\hline\hline
\mu&c_0&c_1&c_0&c_1&\ldots&c_0&c_1&c_0\\\hline
\mu'&c_1&c_0&c_1&c_0&\ldots&c_1&c_0&c_1\\\hline\hline
\end{array}
\end{eqnarray}
If $k$ is even,  the display above must be replaced by the following one:
\begin{eqnarray}\label{joder2}\begin{array}{||c||c|c|c|c|c|c|c|c||}\hline\hline
   G&v_0&e_1&e_2&e_3&\ldots&e_{k-1}&e_k&v_k\\\hline\hline 
\mu&c_0&c_1&c_0&c_1&\ldots&c_1&c_0&c_1\\\hline
\mu'&c_1&c_0&c_1&c_0&\ldots&c_0&c_1&c_0\\\hline\hline
\end{array}
\end{eqnarray}

\begin{figure}[htp]
\hspace*{1.8cm}
\includegraphics[scale=0.9]{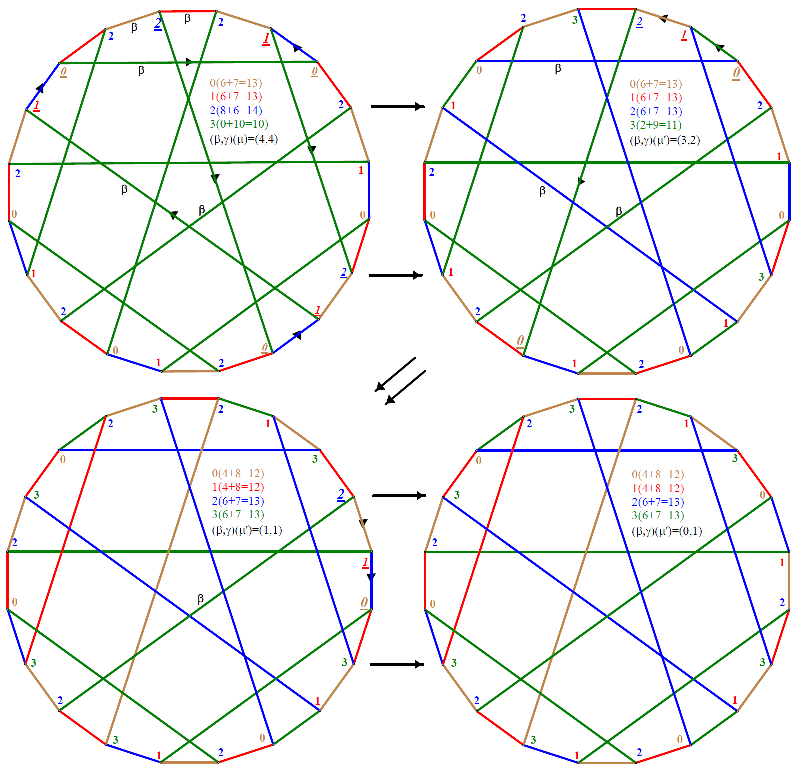}
\caption{From a lacunar sTC to an EqTC for the Desargues graph (Example~\ref{green}).}
\label{figdes}
\end{figure}

In either case, $S=\{v_0,e_1,v_1,e_2,v_2,\ldots,e_k,v_k\}$, is maximal as a color alternating path of $G$ wrt $\mu$. Then, $G$ has all its edges that are not in $S$ but are incident to either $v_0$ or $v_k$ necessarily colored in $\chi(\Delta)\setminus\{c_0,c_1\}$. Thus, the transpositional exchanges between $c_0$ and $c_1$ schematized in displays (\ref{joder1}) and (\ref{joder2}) in passing from $\mu$ to $\mu'$ guarantee that $\mu'$ is a well-defined sTC with $\mu'\neq\mu$,
no matter what colors of $\chi(\Delta)$ the vertices $v_1,\ldots,v_{k-1}$ in $S$ are assigned by $\mu$, and however the possible appearance of new $\beta$-edges in $\mu'$ that were not in $\mu$. like what happens to the edge $vc$ in Example~\ref{joder3}. In any case, the $\mu$-colors of the vertices in $V(G)\setminus\{v_0,v_k\}$ are preserved by $\mu'$, which differs from $\mu$ only on the elements of the skeleton $\rho(S)=\{v_0.e_1,v_1,\ldots,e_k,v_k\}$. 
It is clear that if the set of $\beta$-edges of $\mu$ strictly contains the set of $\beta$-edges of $\mu'$, then $\beta(\mu')<\beta(\mu)$, which establishes item (a). Finally, it is equivalent to say that the difference between the maximum and minimum cardinalities is less for the $\mu'$-classes than for the $\mu$-classes that to write $\gamma(\mu')<\gamma(\mu)$; this establishes item (b).
\end{proof}

Given a cycle $\eta$, it is easy to construct an sTC $\mu_\eta$ on $\eta$ via color set $\{0,1,2\}=\chi(2)\subset\chi(3)$.
 The following definition of {\it lacunar} TC or sTC $\mu$ of a connected simple graph $G$ includes the case of $G$ cubic having a Hamilton cycle $\eta$. 
 Any sTC $\mu_\eta$ on such $\eta$ determines an sTC $\mu$ on $G$ as an extension of $\mu_\eta$ so that the remaining 1-factor in $G$ after removing $\eta$ can have its edges with color $3\in\chi(3)$ forming part of $\mu$. Theorem~\ref{t1} is applicable to such situations, which is a frequent case in the applications ahead in Sections~\ref{s2}-\ref{last}-\ref{s3}.

\begin{definition}
A TC or  sTC $\mu:V(G)\cup E(G)\rightarrow\chi(\Delta)$ of a connected simple graph $G$ with largest degree $\Delta$ is said to be {\it lacunar} if there exists a color $c\in\chi(\Delta)$ such that no $v\in V(G)$ has $\mu(v)=c$.
\end{definition}

\begin{definition}\label{d1} Let $G$, $\mu$, $S$ and $\mu'$ be as in the statement of Theorem~\ref{t1}. Then,  $\mu'$ is said to be a {\it partial $\beta$-reduction}
of $\mu$ wrt $S$ if $0<\beta(\mu')$, and  
it is said to be a {\it total $\beta$-reduction} of $\mu$ wrt $S$ if $0=\beta(\mu')$. 
On the other hand, $\mu'$ is said to be a partial {\it $\gamma$-reduction} of $\mu$ wrt $S$, if $1<\gamma(\mu')$, and it is said to be a {\it total $\gamma$-reduction} of $\mu$ wrt $S$ if either $0=\gamma(\mu')$, in case $\mu'$ is a TC, or $1=\gamma(\mu')$, otherwise.
\end{definition}

\begin{figure}[htp]
\includegraphics[scale=0.56]{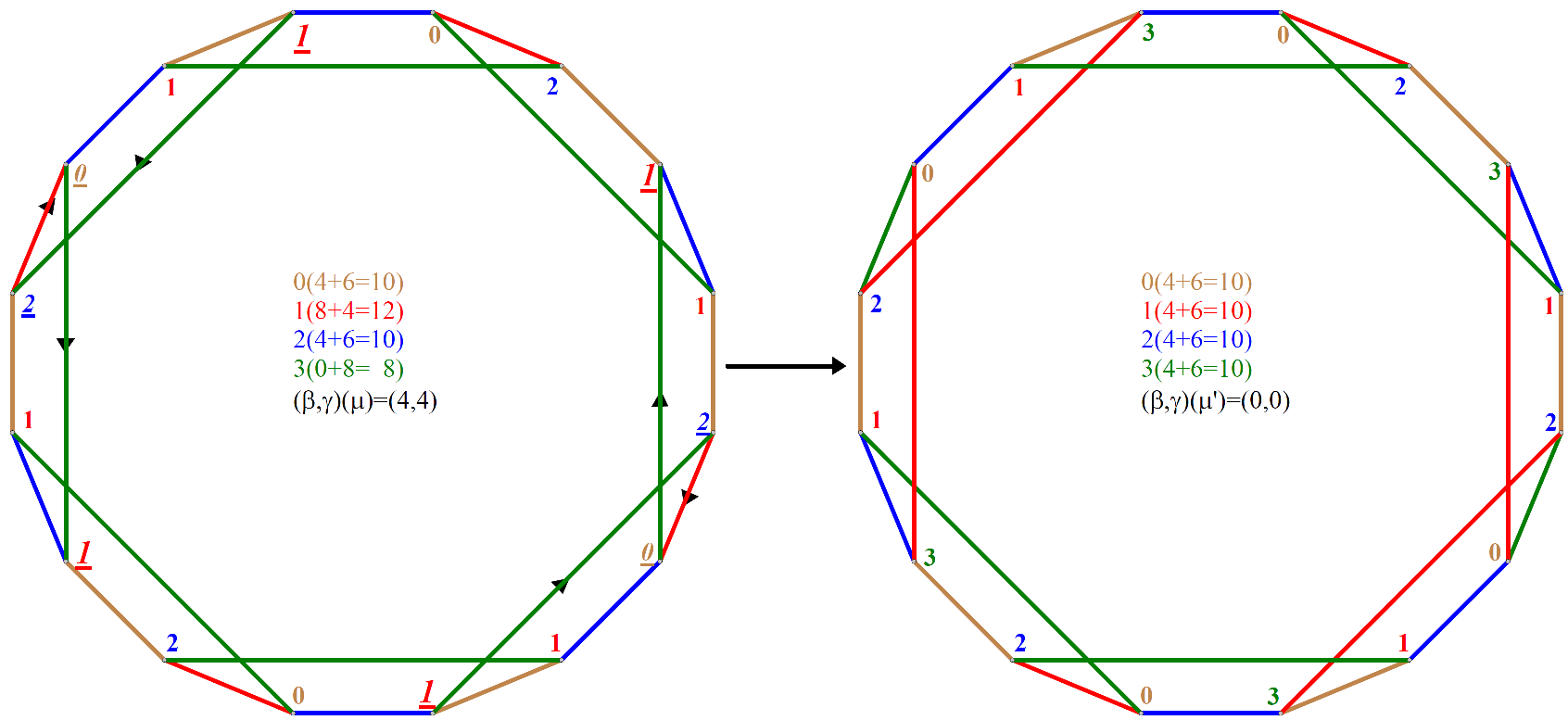}
\caption{From a lacunar sTC to an EqTC for $C_8\square K_2$ (Example~\ref{c8k2}).}
\label{tai}
\end{figure}

\begin{definition}\label{d2}
Let $M=(\mu_0,\mu_1,\ldots,\mu_t)$ be a sequence of sTCs ($0<t\in\mathbb{Z}$) such that $\mu_i$ is a partial $\delta_i$-reduction of $\mu_{i-1}$ wrt the MCAP $S_{i-1}$ of $G$, for $0<i<t$, and $\mu_t$ is a total $\delta_t$-reduction of $\mu_{t-1}$ wrt the MCAP $S_{t-1}$ of $G$, where either $\delta_i=\beta$ or $\delta_i=\gamma$, for $0<i\le t$. Then, $M$ is said to be an {\it sTC-to-EqTC reduction sequence} of $\mu_0$ wrt $(S_0,S_1,\ldots,S_{t-1})$ if all values of 
$\delta_1,\ldots,\delta_t$ are either $\beta$ or $\gamma$. A finite sequence of $\beta$-reductions followed by a finite sequence of $\gamma$-reductions is said to be a {\it $\beta$-$\gamma$-reduction sequence}. 
\end{definition}

\subsection{TCs via covering graph maps}\label{cov}


\begin{definition}
Let $G=(V,E)$ be a graph and let $v\in V$. The {\it open neighborhood} $N_G(v)$ of $v$ in $G$ is the set of all vertices of $G$ adjacent to $v$. Let $G'=(V',E')$ be a graph and let $f=(f_V,f_E):G'\rightarrow G$ be a graph map composed by surjections $f_V:V'\rightarrow V$ and $f_E:E'\rightarrow E$ such that $f_E(v'w')=f_V(v')f_V(w'), \forall v'w'\in E'$. If $f|_{N_{G'}(v')}$ is bijective onto $N_{G}(f(v')), \forall v'\in V'$, then $G'$ is said to be a {\it cover} of $G$ and $f$ a {\it covering graph map} of $G$. If $\exists\,  0<r\in\mathbb{Z}\backepsilon |f^{-1}(v)|=r, \forall v\in V$, then $f$ is said to be {\it $r$-fold} and $G'$  an {\it $r$-cover}. If $\Delta(G')=\Delta(G)=\Delta$ and $\mu:G\rightarrow\chi(\Delta)$ is an sTC, then an sTC $\mu':G'\rightarrow\chi(\Delta)$ is said to be a {\it lifting} of $\mu$ if 
$\mu'(v')=\mu(f_V(v'))$, $\forall v'\in V'$, and $\mu'(e')=\mu(f_E(e')$, $\forall e'\in E'$.
\end{definition}

\begin{figure}[htp]
\includegraphics[scale=0.85]{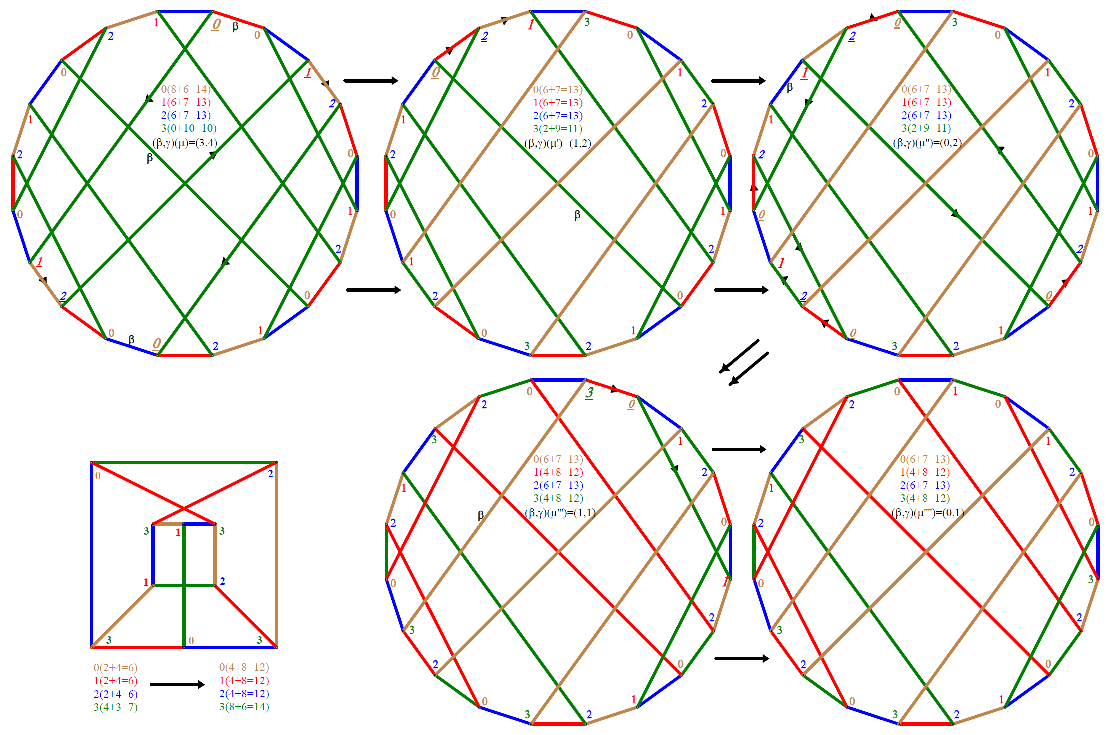}
\caption{EqTC in $Pet$ and $Dod$ passing from lacunar sTC to EqTC (Example~\ref{pete}).}
\label{petdod}
\end{figure}

\begin{theorem}\label{fold}
Given an $r$-fold covering graph map $f:G'\rightarrow G$ of connected simple graphs $G'$ and $G$, both with maximum degree $\Delta$, and given an sTC $\mu:G\rightarrow\chi(\Delta)$, then there exists a corresponding sTC $\mu':G'\rightarrow\chi(\Delta)$ that is a lifting of $\mu$
such that $\beta(\mu')=r\beta(\mu)$ and $\gamma(\mu')=r\gamma(\mu)$. If $\mu$ is a TC, so is $\mu'$. If $\mu$ is an EqTC and $\gamma(\mu)=0$, then $\mu'$ is an EqTC.   
\end{theorem}

\begin{proof}
Consider the inverse image $v'$ or $e'$ of any vertex $v$ or edge $e$ of $G$ and assign the color $\mu(v)$ or $\mu(e)$ to $v'$ or $e'$, respectively, expressed as $\mu'(v')=\mu(v)$ or $\mu'(e')=\mu(e)$. In addition, each $\beta$-edge $e$ of $G$ has $|f^{-1}(e)|=r$. Moreover, if $\mu$ is equitable and $\gamma(\mu)=1$, then $\gamma(\mu')=r$, so that if $r>1$, then $\mu'$ is not equitable; if $\gamma(\mu)=0$, then $\gamma(\mu')=0$ and $\mu'$ is equitable.\end{proof}

\subsection{Perfect codes and total-perfect codes}\label{perfect}

\begin{definition}\label{totalp} Let $G=(V,E)$ be a connected simple graph. 
An {\it efficient dominating set}, or {\it perfect code} \cite{D73} of $G$ is an independent subset $S$ of $V$ such that $|N(v)\cap S|=1$, $\forall v\in V\setminus S$. 
A {\it total-perfect code} \cite{KW} of $G$ is a subset $S$ of $V$ such that $|N(v)\cap S|=1$, $\forall v\in V$. 
If $\mu$ is a TC for which the vertex sets of the $\mu$-classes are perfect codes, then $\mu$ is said to be an {\it efficient TC} \cite{prev}.
If $G$ is cubic with the vertex sets of its $\mu$-classes being 
total-perfect codes, then $\mu$ is said to be {\it 3-total-perfect} if each such code is {\it imperfect}, meaning it induces a matching of $\beta$-edges, and is said to be {\it 1-total-perfect} if only one of the three vertex sets of the $\mu$-classes is an imperfect total-perfect code.
 \end{definition}

\begin{definition}
Let $2<r\in\mathbb{Z}$. The {\it M\"obius ladder} $Mob_r$ is the graph formed from the $2r$-cycle $C_{2r}$, called its {\it external cycle}, by adding edges connecting its opposite vertices.
\end{definition}

\begin{definition}\cite{Coxc}
A {\it cyclic configuration} $n_k$ consists of $n$ points $0,1,\ldots,n-1$  and $n$ lines, each composed by $k $ points, each point as the intersection of $k$ lines, subject to symmetry mod $n$. The {\it Levi graph} of $n_k$ is the graph whose vertices are the points and lines of $n_k$ with an edge $p\ell$ whenever the point $p$ is in the line $\ell$.
\end{definition} 

\begin{definition}\label{n_k} \cite{H,PiRa} Let $0<N=b_02^{n-1}+\cdots+b_{n-2}2+b_{n-1}\in\mathbb{Z}$, so $N$ is the binary number $b_0b_1\cdots b_{n-1}$. The {\it Haar graph} $H(N)$ is the bipartite 
with ''white" vertex part $\{v_0,\ldots,v_{n-1}\}$ and ''black" vertex part $\{u_0,\ldots,u_{n-1}\}$, and an edge $u_{i+j}v_i$, (indices taken mod $n$), for each $i=0,1,\ldots,n-1$ and each $j$ for which $b_j=1$, ($j=0,1,\ldots,n-1$).
\end{definition}

\begin{figure}[htp]
\includegraphics[scale=0.84]{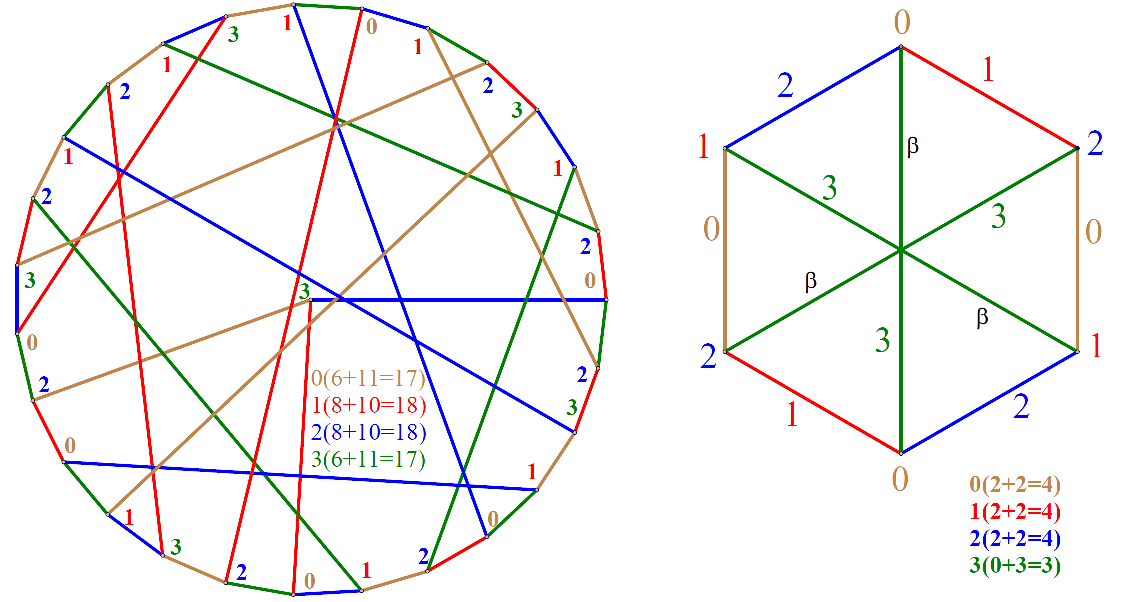}
\caption{EqTC of $Cox$ (Example~\ref{coxa}) and lacunar equitable sTC of $K_{3,3}$ (Example~\ref{K33}).}
\label{coxk33}
\end{figure}

A cyclic configuration $n_k$ is determined by a specific line $B\subseteq\{0,1,\ldots,n-1\}$ of $n_k$, where $|B|=k$. Such $B$ is in 1-1 correspondence with the integer $N=b_02^{n-1}+\cdots+b_{n-2}2+b_{n-1}$, where $b_i=1$ if $i\in B$ and $b_i=0$ if $i\notin B$. 
In Definition~\ref{n_k}, $H(N)$ is the Levi graph of the cyclic configuration $n_k$. Moreover,
Levi graphs of cyclic configurations are precisely the Haar graphs that have girth at least 6 \cite{Boben,H}.

\begin{definition}
Let $2<r\in\mathbb{Z}$. Let $C_{2r}=(v_0,e_1,v_1,e_2\ldots,e_r,v_r,\ldots,e_0)$ be the external cycle of $Mob_r$. Let us replace the total 2-paths $v_{i-1},e_i,v_i,e_{i+1},v_{i+1}$ and $v_{r+i-1},e_{r+i},v_{r+i},e_{r+i+1},\\ v_{r+i+1}$ and the edge $v_iv_{r+i}$, where $i=0,1,\ldots,r-1\mod 2r$, by the total 3-paths $v_{i-1},e_i,v'_i,\\ e'_i,v''_{i+1},e_{i+1},v_{i+1}$ and  $v_{r+i-1},e_{r+i},v'_{r+i}, e'_{r+i},v''_{r+i+1},e_{r+i+1},v_{r+i+1}$ and edge pair $v'_iv''_{r+i}$ and $v''_iv'_{r+i}$, respectively. The resulting graph is the Haar graph $H(2^{r-2}(2^{r-1}+1)+1)$, (where decimal $2^{r-2}(2^{r-1}+1)+1$ is binary $(10^{r-3}10^{r-4}1)_2$),
 and will be said to be the {\it Fat-M\"obius ladder}, denoted $FMob_r$.
\end{definition}
 
\section{Applications to symmetric cubic graphs}\label{s2}

Here are some applications of Theorem~\ref{t1} and Definitions~\ref{d1} and~\ref{d2} to symmetric cubic graphs, meaning graphs that are both vertex-transitive and edge-transitive.
Cases of cubic hamiltonian graphs may use the LCF notation for their expression \cite{CFP,Frucht}. In the cases of distance transitive graphs, their intersection arrays \cite{brower} will be mentioned, too.

\subsection{Total reductions into EqTCs}

\begin{figure}[htp]
\includegraphics[scale=0.8]{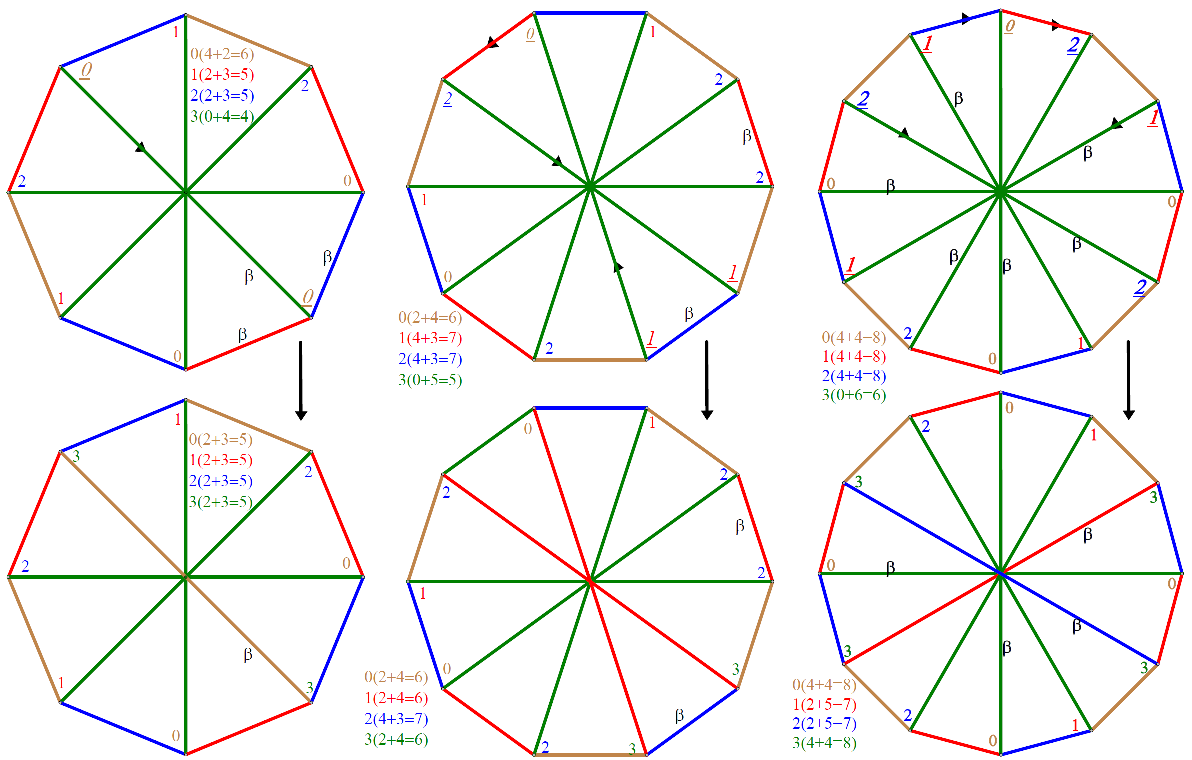}
\caption{From lacunar sTCs to equitable sTCs of $Mob_r$ ($r=4,5,6$). (Example~\ref{K33}).}
\label{mob12}
\end{figure}
 
\begin{example}\label{q3}
The distance transitive 3-cube graph $Q_3$ with intersection array $\{3,2,1;1,2,3\}$ is the hamiltonian Haar graph $H(11)=H(1011_2)$ with LCF notation $[-3,3]^4$, seen in Figure~\ref{cubo} as the union of an external Hamilton 8-cycle $\eta$ and an internal 1-factor $F$ formed by 4 edges. $Q_3$ is represented on the left of Figure~\ref{cubo} with a lacunar sTC $\mu$ that uses color set $\chi(\Delta)=\chi(3)=\{0=$ hazel, $1=$ red, $2=$ blue, $3=$ green$\}$. This way, $\eta$ is seen clockwise from the upper-left vertex via the color cycle
$\mu(\eta)=(1_20_12_00_12_01_21_02_10_21_0)$, where edge colors appear as subindices of their preceding vertex colors, while $F$ has only green (color 3) edges.
There are two $\beta$-edges for $\mu$, marked with the Greek letter $\beta$. 

\begin{figure}[htp]
\includegraphics[scale=0.85]{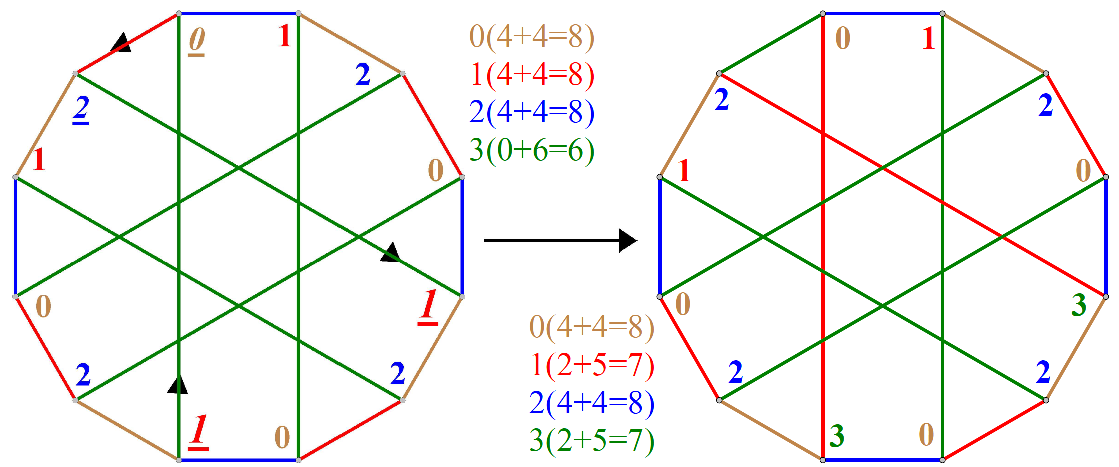}
\caption{From a lacunar TC to an EqTC for the Franklin graph $H(11)$ (Example~\ref{fr}).}
\label{frank}
\end{figure}

With initial vertex $w_0$ being the one corresponding to the upper-left vertex color red (1) in $\eta$,
an MCAP $S=(w_0,e_1,w_1,e_2,w_2,e_3,w_3)$ with color path $\mu(S)=1_32_10_31$ is distinguished by having its vertex color numbers in underlined Italics in Figure~\ref{cubo} (in contrast with the remaining vertex numbers, in Roman type), its two degree-one vertices in larger type and its edges oriented in the left-to-right direction of $S$, so that $\mu(\rho(S))=\mu(w_0,e_1,e_2,e_3,w_3)=(1,3,1,3,1)$. The right side of Figure~\ref{figshp} shows the TC $\mu'$ for $Q_3$ provided by Theorem~\ref{t1}. It is obtained by the ongoing $\beta$-reduction, a total $\beta$-reduction. We note that $\mu(S)=1_32_10_31$ becomes on this right side  into $\mu'(S)=3_12_30_13$. In the center of such right side, a listing showing the cardinalities of the $\mu'$-classes shows that this TC is in fact an EqTC, where each row of the listing is of the form $c(|V^c|+|E^c|=|Q_3^c|)$, with $c\in\chi(3)$ and $|V^c|,|E^c|$ and $|Q_3^c|$ standing for the cardinalities of the vertex, edge and element (vertex or edge) subsets of the $\mu'$-class of $c$. This listing thus is of the form $i(2+3=5)$, for $i\in\chi(\Delta)$, with common cardinality 5. We resume this by writing $Q_3(|Q_3^0|,|Q_3^1|,|Q_3^2|,|Q_3^3|)=Q_3(5,5,5,5)=Q_3(5^4)$. As a result, no $\gamma$-reduction for $\mu'$ is necessary here. In sum, Figure~\ref{cubo} shows a $\beta$-$\gamma$-reduction sequence from a lacunar sTC $\mu$ onto an EqTC $\mu'$ of $Q_3$, passing from $(\beta,\gamma)(\mu)=(\beta(\mu),\gamma(\mu))=(2,2)$ to $(\beta,\gamma)(\mu')=(0,0)$. 
\end{example}

\begin{remark}
The right side of Figure~\ref{cubo} is  represented as a cutout on the left of the symbol ``$\neq$" in display (\ref{oct}), where vertical color 1 edges must be identified. Note that
pairs of antipodal vertices in $Q_3$ (at distance 3) are assigned the same color number in the TCs on either side of ``$\neq$", (where vertical color 3 edges must be identified on the cutout of $Q_3$ on the right of ``$\neq$"), so that these two TCs are efficient TCs. Vertex colors are similar on both cutouts but the edge colors differ, so the two edge colorings are {\it orthogonal} \cite[Def. 5]{prev}.  
 
\begin{eqnarray}\label{oct}\begin{array}{lll}
0\hspace*{2.2mm} _-^3\hspace*{2.2mm} 1\hspace*{2.2mm} _-^0\hspace*{2.2mm} 2\hspace*{2.2mm} _-^1\hspace*{2.2mm} 3\hspace*{2.2mm} _-^2\hspace*{2.2mm} 0&&
0\hspace*{2.2mm} _-^2\hspace*{2.2mm} 1\hspace*{2.2mm} _-^3\hspace*{2.2mm} 2\hspace*{2.2mm} _-^0\hspace*{2.2mm} 3\hspace*{2.2mm} _-^1\hspace*{2.2mm} 0\\

\!_1|\hspace*{6mm}_2|\hspace*{6mm}_3|\hspace*{6mm}_0|\hspace*{6mm}_1|&\neq&
\!_3|\hspace*{6mm}_0|\hspace*{6mm}_1|\hspace*{6mm}_2|\hspace*{6mm}_3|\\

2\hspace*{2.2mm} _0^-\hspace*{2.2mm}3\hspace*{2.2mm} _1^-\hspace*{2.2mm} 0\hspace*{2.2mm} _2^-\hspace*{2.2mm} 1\hspace*{2.2mm} _3^-\hspace*{2.2mm} 2&&
2\hspace*{2.2mm} _1^-\hspace*{2.2mm}3\hspace*{2.2mm} _2^-\hspace*{2.2mm} 0\hspace*{2.2mm} _3^-\hspace*{2.2mm} 1\hspace*{2.2mm} _0^-\hspace*{2.2mm} 2 \\
\end{array}\end{eqnarray}
\end{remark}

\begin{theorem}\label{girth}(\cite[Theorem 12]{prev})
The 3-cube $Q_3$ and the prisms $C_{4j}\square K_2$ ($1<j\in\mathbb{Z}$) have efficient TCs in which every 4-cycle uses the four colors of $\chi(3)$ both on vertices and also on edges in two edge mutually orthogonal ways. Moreover, any two such orthogonal ways guarantee an edge coloring in the corresponding 4-cube $Q_4=Q_3\square K_2$, or prism $(C_{4j}\square K_2)\square K_2$ in which the edges of each 4-cycle uses the 4 colors of $\chi(3)$. 
\end{theorem}

\begin{proof} The proof of the examples in this result was already presented in \cite[Theorem 12]{prev}. 
\end{proof}

\begin{example}\label{e1} The distance transitive Heawood graph $Hea=H(69)$ with intersection array $\{3,2,2;1,1,3\}$ is hamiltonian with LCF notation $[5,-5]^7$ seen (twice) in Figure~\ref{figshp} as the union of an external Hamilton 14-cycle $\eta$ and an internal 1-factor $F$ of 7 edges. $Hea$ is represented on the left side of Figure~\ref{figshp} with a lacunar sTC $\mu$ that uses color set $\chi(\Delta)=\chi(3)$. This way, $\eta$ is seen clockwise from the upper-left vertex via the color cycle
$\mu(\eta)=(1_02_10_2(0_12_01_2)^30_10_2)$, where edge colors appear as subindices of their preceding vertex colors, and $F$ has only green (color 3) edges.
There are two $\beta$-edges for $\mu$, distinguished via Greek letter $\beta$. 

\begin{figure}[htp]
\includegraphics[scale=0.86]{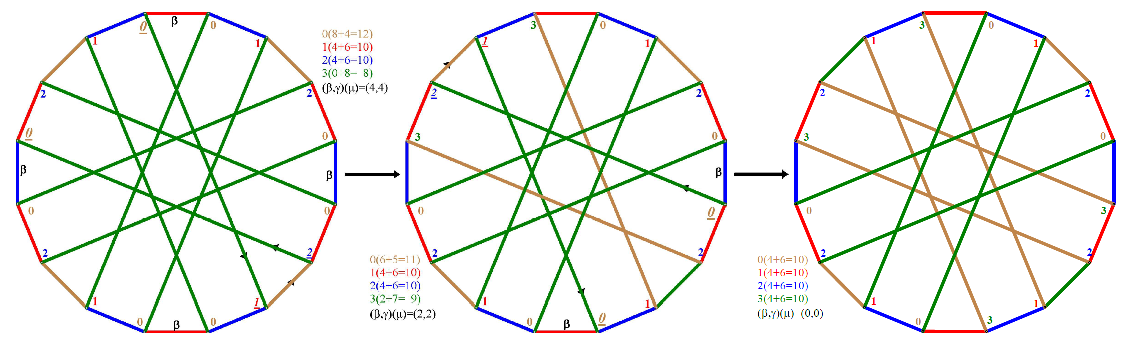}
\caption{Fat-M\"obius graph $FMob_4=H(37)$ (Example~\ref{fatty}).}
\label{fat}
\end{figure}

With initial vertex $w_0$ being the one corresponding to the third vertex color hazel (0) in $\eta$,
an MCAP $S=(w_0,e_1,w_1,e_2,w_2,e_3,w_3,e_4,w_4,e_5,w_5)$ with color path $\mu(S)=0_31_02_31_02_30$ is distinguished, too, by having its vertex color numbers in Italics in Figure~\ref{figshp} (in contrast with the remaining vertex numbers, in Roman type), its two degree-one vertices in larger type and its edges oriented in the left-to-right direction of $S$, so that $\mu(\rho(S))=\mu(w_0,e_1,e_2,e_3,e_4,e_5,w_5)=(0,3,0,3,0,3,0)$. The right side of Figure~\ref{figshp} shows the  TC $\mu'$ for $Hea$ provided by Theorem~\ref{t1}. It was obtained by the ongoing $\beta$-reduction, a total $\beta$-reduction. We note that $\mu(S)=0_31_02_31_02_30$ becomes on the right side of Figure~\ref{figshp} into $\mu'(S)=3_01_32_01_32_03$. In the center of this right side of Figure~\ref{figshp}, a listing showing the cardinalities of the $\mu'$-classes shows that this TC is in fact an EqTC, where each row of the listing is of the form $c(|V^c|+|E^c|=|Hea^c|)$, with $c\in\chi(3)$ standing for color and $|V^c|,|E^c|$ and $|Hea^c|$ standing for the cardinalities of the vertex, edge and element (vertex or edge) subsets of the $\mu'$-class of $c$. The listing thus is $i(4+5=9)$, for $i=0,1,2$ and $3(2+6=8)$, with cardinalities differing at most in one unit. 
 We resume this by writing $Hea(|Hea^0|,|Hea^1|,|Hea^2|,|Hea^3|)=Hea(9,9,9,8)=Hea(89^3)$.As a result, no need of a $\gamma$-reduction for $\mu'$ is necessary here.
\end{example}

\begin{figure}[htp]
\hspace*{1.4cm}
\includegraphics[scale=0.7]{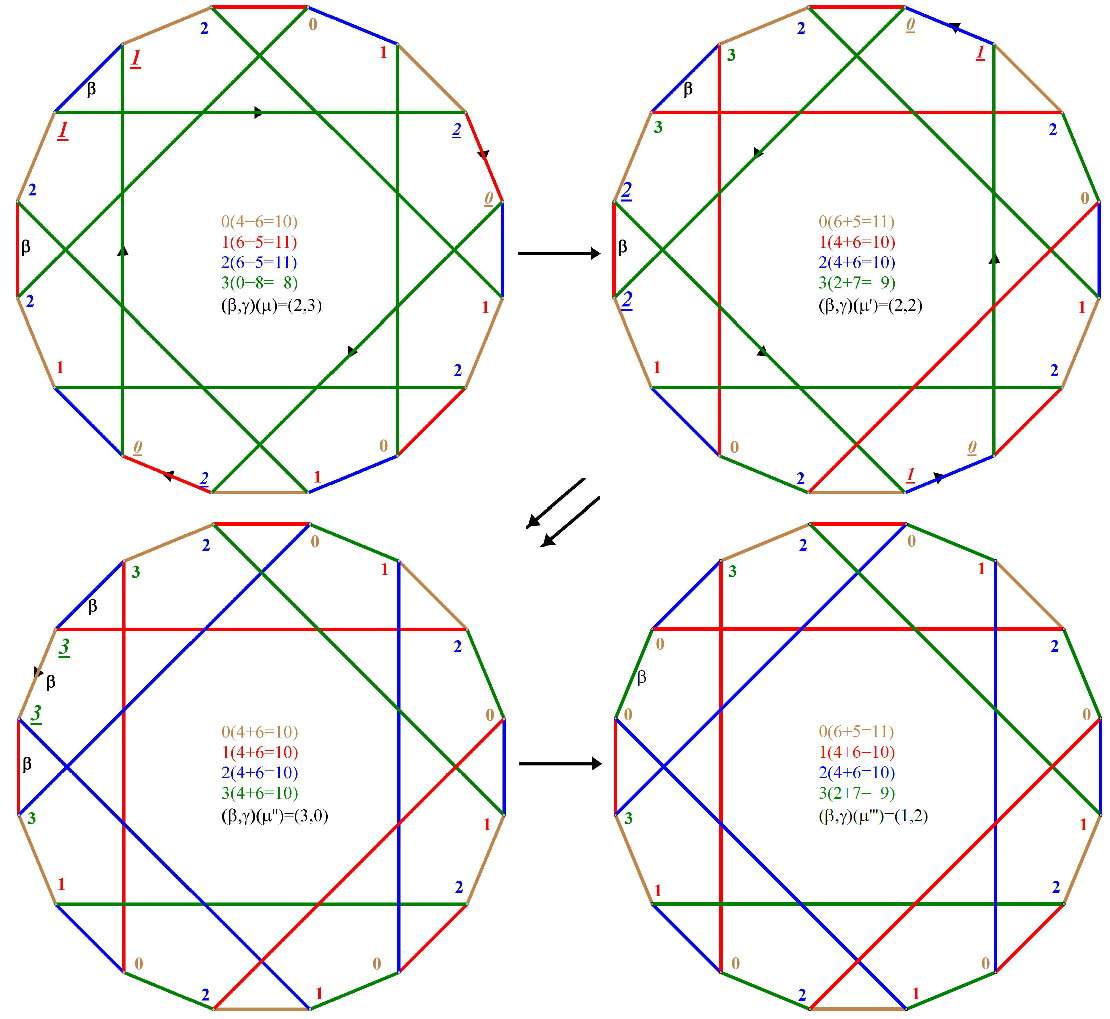}
\caption{M\"obius-Kantor graph $H(133)$ (Example~\ref{moka}).}
\label{moka-joe}
\end{figure} 

\begin{example}\label{Pap} 
The distance transitive Pappus graph $Pap$ with intersection array $\{3,2,2,1;$  $1,2,3\}$ is hamiltonian with LCF notation $[5,7,-7,7,-7,-5]^3$ seen (four times) in Figure~\ref{figeqp} as the union of an external Hamilton 18-cycle $\eta$ and an internal 1-factor $F$ of 9 edges. $Pap$ is represented on the upper-left quarter of Figure~\ref{figeqp} with a lacunar sTC $\mu$ that uses color set $\chi(3)$. There, $\eta$ is seen clockwise from the upper-left vertex via the color cycle
$\mu(\eta)=(1_21_02_10_2(0_12_01_2)^40_22_0)$, where edge colors appear as subindices of their preceding vertex colors, the exponent 4 means concatenation of $0_12_01_2$ four times and $F$ has only green (color 3) edges.

With initial vertices $w_0$ being the ones corresponding to the upper-right, leftmost and again upper-right vertices, colored red (1), blue (2) and hazel (0), respectively,  in $\eta$, successive
MCAPs $S_0=(w_0,e_1,w_1,e_2,w_2)$ and $S_1=(w_0,e_1,w_1e_2,w_2,e_3,w_3)$ and $S_2=(w_0,e_1,e_2,e_3)$ with color paths $\mu(S_0)=1_02_10$, $\mu'(S_1)=2_30_21_32$ and $\mu''(S_1)=0_32_01_30$ are distinguished, too, by having their vertex color numbers in underlined Italics in the upper-left, upper-right and lower-left quarters of Figure~\ref{figeqp} (in contrast with the remaining vertex numbers, in Roman type), its degree-one vertex pairs in larger type and their edges oriented in the direction of the three respective MCAPs, so that $\mu(\rho(S_0))=1010$, $\mu'(\rho(S_1))=(2,3,2,3,2)$, $\mu''(\rho(S_2))=(0,3,0,3,0 )$. 

The lower-right quarter of Figure~\ref{figeqp} shows the EqTC $\mu'''$ of $Pap$ provided by Theorem~\ref{t1}. It was obtained by an initial $\beta$-reduction and two subsequent $\gamma$-reductions, composing an sTC-to-EqTC reduction sequence departing as said from a lacunar sTC and passing successively through a lacunar TC with automorphism group $\mathbb{Z}_6$ and one more TC. In the center of each of the last three quarters in Figure~\ref{figeqp}, a listing showing the cardinalities of the $\mu'$-classes shows that this TC is in fact an EqTC, in a fashion similar to the listing on right side of Figure~\ref{figshp}. The final listing is $i(4+7=11)$, for $i=0,2,3$ and $1(6+6=12)$, with cardinalities differing at most in one unit.  We resume this by writing $Pap(|Pap^0|,|Pap^1|,|Pap^2|,|Pap^3|)=Pap(11,11,11,12)=Pap(11^312)$.
\end{example}

\begin{figure}[htp]
\includegraphics[scale=0.89]{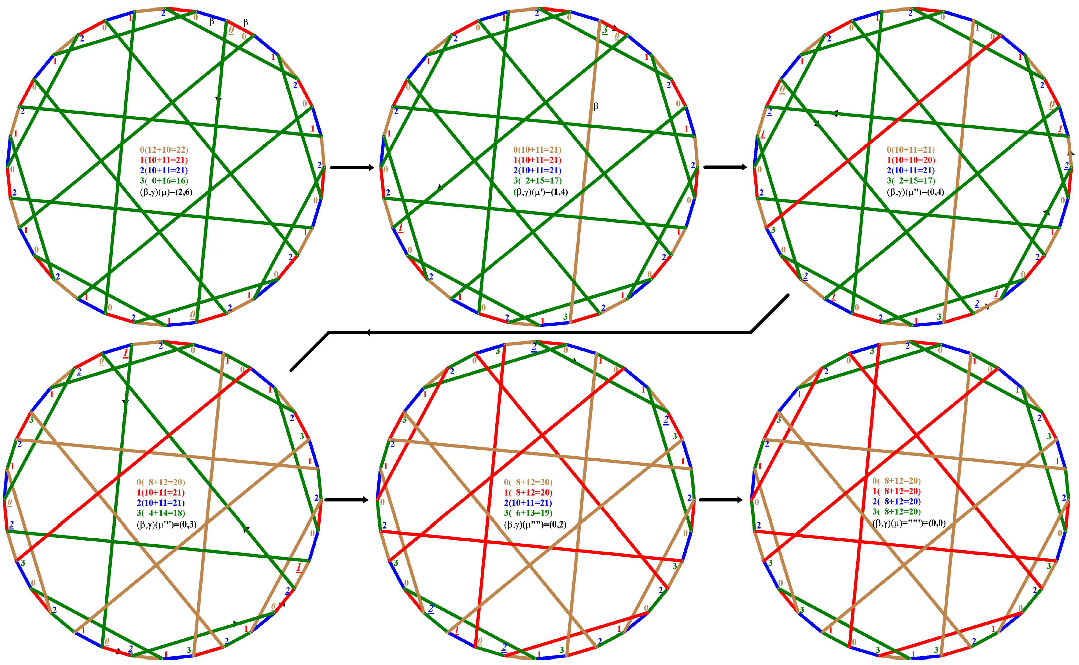}
\caption{From an sTC of the Dyck 32-vertex graph $Dy$ to an EqTC (Example~\ref{dy}).}
\label{figdyck}
\end{figure}

\begin{example}\label{green} 
The distance transitive Desargues graph $Des$ with intersection array $\{3,2,2,$ $1,1,1;1,1,2,2,3\}$ is hamiltonian with LCF notation $[5,-5,9,-9]^5$ seen (four times) in Figure~\ref{figdes} as the union of an external Hamilton 20-cycle $\eta$ and an internal 1-factor $F$ of 10 edges. $Des$ is represented on the upper-left quarter of Figure~\ref{figdes}  with a lacunar sTC $\mu$ that uses color set $\chi(3)$. There, $\eta$ is seen clockwise from the upper-left vertex via the color cycle $\mu(\eta)=(2_1(2_01_20_1)^62_0)$, where edge colors appear as subindices of their preceding vertex colors and $F$ has only green (color 3) edges.

With initial vertices $w_0$ being the ones corresponding to the upper-left (or first), fourth and fifth vertices, colored blue (2), hazel (0) and again blue (2), respectively,  in $\eta$, successive
MCAPs $S_0=(w_0,e_1,w_1,e_2,w_2,e_3,w_3,e_4,w_4,e_5,w_5,e_6,w_6,e_7,w_7)$, $S_1=(w_0,e_1,w_1,e_2,w_2,e_3,$ $w_3)$, $S_2=(w_0,e_1,w_1,e_2,w_2)$ with color paths $\mu(S_0)=2_30_21_31_20_30_21$, $\mu'(S_1)=0_31_02_30$, $\mu''(S_2)=2_01_20$ are distinguished, too, by having their vertex color numbers in underlined Italics in the upper-left, upper-right and lower-left quarters of Figure~\ref{figdes} (in contrast with the remaining vertex numbers, in Roman type), their degree-one vertex pairs in larger type  and their edges oriented in the direction of the three respective MCAPs, so that $\mu(\rho(S_0))=(2,3,2,3,2,3,2,3,2)$, $\mu'(\rho(S_1))=(0,3,0,3,0 )$, $\mu''(\rho(S_2))=(2,0,2,0)$. The lower-right quarter of Figure~\ref{figdes} shows a EqTC $\mu'''$ of $Des$ provided by Theorem~\ref{t1}. It was obtained by three successive $\beta$-reductions, the first two of which being $\gamma$-reductions as well, and composing an sTC-to-EqTC reduction sequence. Inside the lower-right quarter of Figure~\ref{figdes}, a listing shows $i(4+8=12)$, for $i=0,1$ and $j(6+7=13)$, for $j=2,3$, with cardinalities differing at most in one unit. 
 We resume this by writing $Des(|Des^0|,|Des^1|,|Des^2|,|Q_3^3|)=Des(12,12,13,13)=Des(12^2,13^2)$.
No $\gamma$-reduction is necessary here, since an EqTC was finally obtained.
\end{example} 

\subsection{TCs via covering graph maps}

\begin{example}\label{c8k2}
For $r>1$, the prism graph $C_{4r}\square K_2$ is an $r$-covering graph of $Q_3$ via an $r$-fold from $C_{4r}\square K_2$, where $C_{4r}$ is a $(4r)$-cycle, so that Theorem~\ref{fold} applies from the TC of $Q_3$ in Example~\ref{q3} represented on the right side of Figure~\ref{cubo}. Figure~\ref{tai} shows a $\beta$-$\gamma$-reduction sequence from a lacunar sTC $\mu$ onto an EqTC $\mu'$ of $C_8\square K_2$, passing from $(\beta,\gamma)(\mu)=(\beta(\mu),\gamma(\mu))=(4,4)$ to $(\beta,\gamma)(\mu')=(0,0)$. 
We notice that $\mu$ on the right of Figure~\ref{cubo} is efficient, result extendable to every prism $C_{4r}\square K_2$, that is: an efficient TC exists on any such prism.  
A similar immediate equitable conclusion like the one for $\mu$ cannot be carried out by departing from the dodecahedral graph as a 2-cover of the Petersen graph, as seen in the following example.
\end{example}

\begin{figure}[htp]
\includegraphics[scale=0.71]{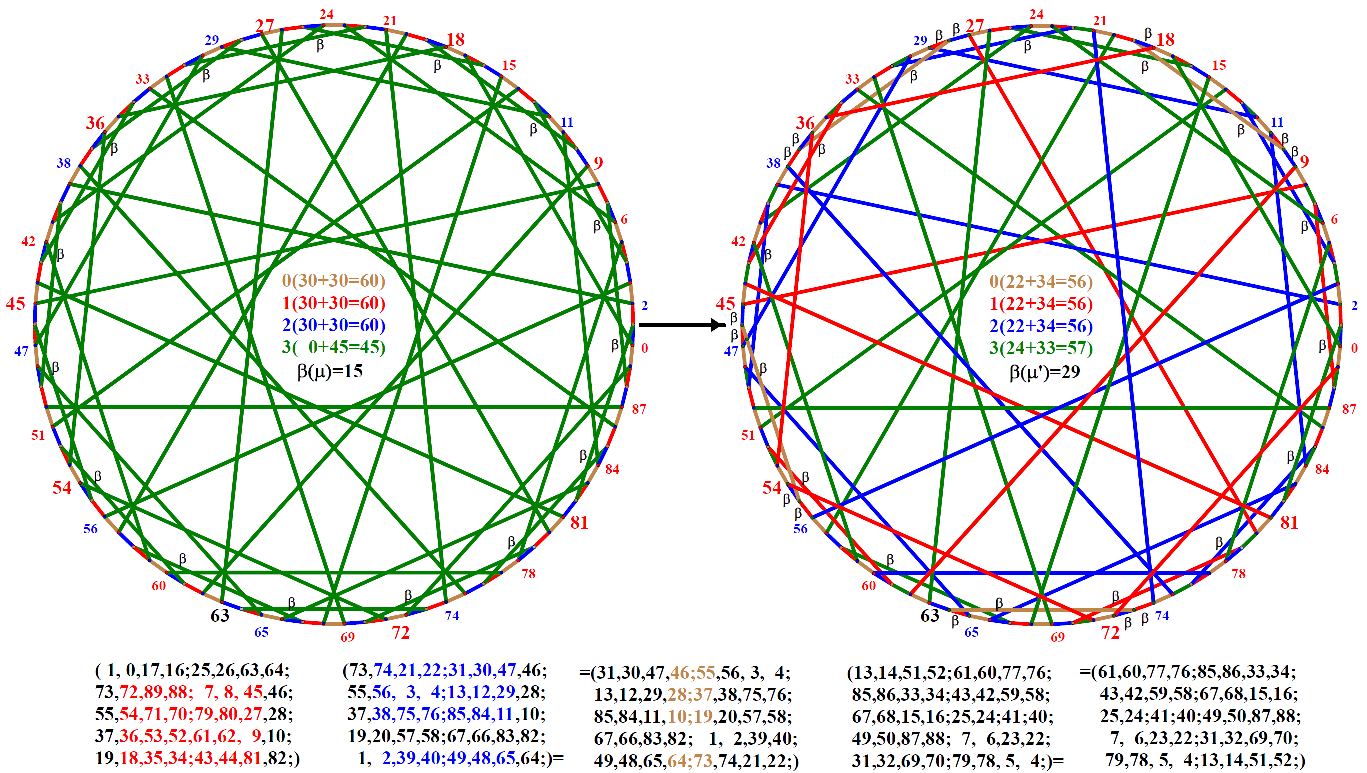}
\caption{From lacunar sTC onto equitable sTC for  Foster 90-vertex graph (Example~\ref{cansado}).}
\label{fig90}
\end{figure}

\begin{example}\label{pete} The Petersen graph $Pet$ is the unique $(3,5)$-cage \cite[p. 175]{Harary} as well as the unique $(3,5)$-Moore graph. The intersection array of $Pet$ is $\{3,2;1,1\}$. It is the only distance transitive cubic graph which is a snark. We use its representation in~\cite[Fig. 1]{O} to show a TC $\mu$ on it in the lower-left of Figure~\ref{petdod}.
 A listing showing the cardinalities of the $\mu'$-classes shows that this TC is in fact an EqTC, where each row of the listing is of the form $c(|V^c|+|E^c|=|Q_3^c|)$, with $c\in\chi(3)$ standing for color and $|V^c|,|E^c|$ and $|Pet^c|$ standing for the cardinalities of the vertex, edge and element (vertex or edge) subsets of the $\mu'$-class of $c$. The listing thus is $i(2+4=6)$, for $i=0,1,2$ and $(3(4+3=7)$ for $i=3$. We resume this by writing $Pet(|Pet^0|,|Pet^1|,|Pet^2|,|Pet^3|)=Pet(6,6,6,7)=Pet(6^37)$. As a result, no  $\gamma$-reduction for $\mu$ is necessary, too. 

The dodecahedral graph $Dod$ is a cubic symmetric graph, namely the generalized Petersen graph $GP(10,2)$ of LCF notation $[10,7,4,-4,-7,10,-4,7,-7,4]^2$ and intersection array $\{3,2,1,1,1;1,1,1,2,3\}$. Theorem~\ref{fold} insures that $Dod$ admits a TC  $\phi^{-1}(\mu)$ (via a covering graph map $\phi:Dod\rightarrow Pet$), however not equitable, as exemplified to the right of the arrow from the listing for $\mu$ in $Pet$ (below its representation in Figure~\ref{petdod}), showing the corresponding listing for $Dod$ and from which is deduced that $\gamma(\phi^{-1}(\mu))=2$. 
 
The rest of Figure~\ref{petdod} illustrates that departing from the lacunar sTC $\mu^0$ of $Dod$ shown on the upper-left of Figure~\ref{petdod} it is feasible to obtain an EqTC $\mu^4$ via a composition $\mu^0\rightarrow\mu^1\rightarrow\mu^2\rightarrow\mu^3\rightarrow\mu^4$
 of two $\beta$-reductions, (followed by) a $\gamma$-reduction and a final $\beta$-reduction. In more detail, an MCAP $S^1=(w_0^1,e_1^1,w_1^1,e_2^1,w_2^1,e_3^1,w_3^1,e_4^1,$ $w_4^1,e_5^1,w_5^1)$ departing from the upper-right vertex of the representation of $Dod$ with color path $\mu(S^1)=0_31_02_31_02_30$) takes as indicated by an arrow to the right to the sTC $\mu_1$ in the upper-middle of Figure~\ref{petdod}. In it, an MCAP $S^2=(w_0^2,e_1^2,w_1^2,e_2^2,w_2^2)$ with color path $\rho(S^2)=0_12_01$ takes to the sTC $\mu^2$ on the upper-right of Figure~\ref{petdod}. In it, an MCAP $S^3=(w_0^3,\ldots,w_9^3)$ with color path $1_30_12_30_12_30_12_30_12_31$ takes to the sTC $\mu^3$ on the lower-middle of Figure~\ref{petdod}. In it, an MCAP $S^4=(w_0^4,e_1^4,w_1^4,e_2^4,w_2$ with color path $\rho(S^4)=3_10_21$ takes to the EqTC $\mu^4$ on  the lower-right of Figure~\ref{petdod}, whose listing inside this last representation of $Dod$ resumes as $Dod(|Dod^0|,|Dod^1|,|Dod^2|,|Dod^3|)=Dod(13,12,13,12)=Dod(12^213^2)$. 
\end{example}
    
\begin{figure}[htp]
\includegraphics[scale=0.87]{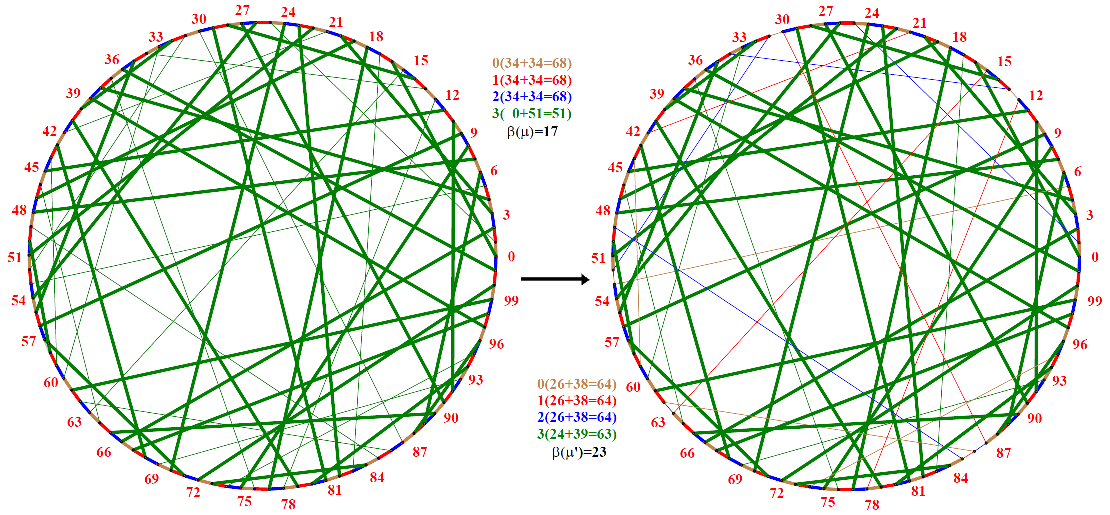}
\caption{Application to the Biggs-Smith graph (Example~\ref{102}).}
\label{fig102}
\end{figure}
  
\begin{example}\label{coxa}
The distance transitive 28-vertex Coxeter graph $Cox$ with intersection array $\{3,2,2,1;1,1,1,2\}$ is shown on the left of Figure~\ref{coxk33} with a TC $\mu$, where a listing is shown, namely $0(6+11=17)$ ,for $i=0,3$, and $i(8+10)=18$, for $i=1,2$.  We resume this by writing $Cox(|Cox^0|,|Cox^1|,|Cox^2|,|Cox^3|)=Cox(17,18,18,17)=Cox(17^218^2)$.
\end{example}

\subsection{Reductions to equitable sTCs and total-perfect codes}\label{equi}

We engage into obtaining from the originally lacunar sTCs treted in this subsection equitable sTCs via $\gamma$-reductions at the possible cost of increasing their beta invariant.

\begin{example}\label{K33}
The distance transitive complete bipartite graph $K_{3,3}=Mob_3=H(7)=H(111_2)$ with intersection array $\{3,2;1,3\}$ has a lacunar sTC $\mu_1$ with $\beta(\mu_1)=3$, represented on the right side of Figure~\ref{coxk33}. 
$K_{3,3}$ is the first case of the M\"obius ladders $Mob_r$ ($r>2$). 
Each $Mob_r$ has a lacunar sTC $\mu_r$. We believe it cannot be $\beta$-reduced to a TC. One can $\gamma$-reduce $\mu_r$ to an equitable (non-lacunar) sTC $\mu'_r$ when $r>2$, as shown in Figure~\ref{mob12} for $r=4,5,6$.  
\end{example}

\begin{definition}
 If $0<r=3k$, where $0<k\in\mathbb{Z}$, and $C_{2r}=(u_0,v_0,w_0,\ldots,u_j,v_j,w_j,\ldots,\\
 u_{2k},v_{2k},w_{2k})$ is a Hamilton cycle of a connected simply cubic graph $G$ with additional 1-factor $F=F_u\cup F_v\cup F_w$ the disjoint union of edge sets $F_z$ given as the products of $k$ transpositions of the sets $\{z_1,\ldots,z_{2k}\}$, where $z\in\{u,v,w\}$, then $G$ is said to be a $6k$-{\it vertex jumble 
 ladder}. In particular, if $F_z=\Pi_{j=1}^k (z_j\,z_{j+k})$, for $z\in\{u,v,w\}$, then $G$ is the {\it M\"obius ladder} $Mob_{3k}$. 
\end{definition}

\begin{theorem}\label{valentino}
Each $6k$-vertex jumble ladder $G$ ($k>0$) contains a lacunar 3-total-perfect sTC $\mu$. 
 Such $\mu$ is equitable if $k=1$, meaning that $K_{3,3}=Mob_3$ is a jumble graph with a lacunar 3-total-perfect sTC $\mu$ that is equitable with $\gamma(\mu)=1$. 
\end{theorem}

\begin{proof} All $6k$-vertex jumble ladders $G$ ($k>0$) admit lacunar sTCs $\mu$ that are 3-total-perfect by having the external Hamilton cycle with color cycle $(0_21_02_1)^{2k}$. Specifically, the induced components of the set of green (color 3) $\beta$-edges form a partition into three total-perfect codes. If $k=1$, then the sTC $\mu$ on $K_{3,3}=Mob_3$ is equitable because
$\gamma(\mu)=1$. In fact, the color cycle pattern $(0_21_02_1)^{2k}$ provides the involved color classes with $2k$ vertices and $2k$ edges, and the remaining (green) class has $3k$ edges, so the color classes yield total sizes $4k,4k,4k,3k$. Hence $\gamma(\mu)=k$. This means that $\gamma(\mu)=1$ if $k=1$.
\end{proof}

\begin{corollary}\label{ladder}
Each M\"obius ladder $Mob_{3k}$ ($k\ge 1$) contains a lacunar 3-total-perfect sTC. 
which is equitable if 
$k=1$. Thus, $K_{3,3}=Mob_3$ is 
a M\"obius ladder with a lacunar 3-total-perfect sTC that is equitable. 
\end{corollary}

\begin{example}
Is $K_{3,3}$ the only lacunar equitable sTC on any connected simple graph? The answer to such question is no, as shown by $P_4$ with consecutive vertex colors $(1,2,2,2)$ and edge colors $(3,1,3)$, providing a lacunar STC with class sizes $(2,3,2)$, so that neither nonempty vertex class is a total perfect code.  
\end{example}

\begin{example}\label{fr}
The Franklin 12-vertex cubic graph $H(11)=FMob_3$ whose embedding on the Klein bottle divides it into regions having a  minimal coloring via six colors yields the sole counterexample to the Heawood Conjecture \cite{Heawood}. $H(11)$ is nonplanar Hamiltonian with LCF notations $[5,-5]^6$ and $[-5,-3,3,5]^3$. 
Figure~\ref{frank} shows a $\gamma$-reduction of $H(11)$ from a lacunar TC into an EqTC.
\end{example}

\begin{example}\label{fatty}
The lacunar sTC $\mu$ of Fat-M\"obius graph $FMob_4=H(37)$ depicted on the left of Figure~\ref{fat} has $(\beta,\gamma){\mu}=(4,4)$ and is $\beta$-$\gamma$-reduced in the middle of Figure~\ref{fat} to an sTC $\mu'$ with $(\beta,\gamma)(\mu')=(2,2)$. This in turn is $\beta$-$\gamma$-reduced on the right of Figure~\ref{fat} to an EqTC $\mu''$ with $(\beta,\gamma)(\mu'')=(0.0)$.
\end{example}

\begin{figure}[htp]
\includegraphics[scale=0.88]{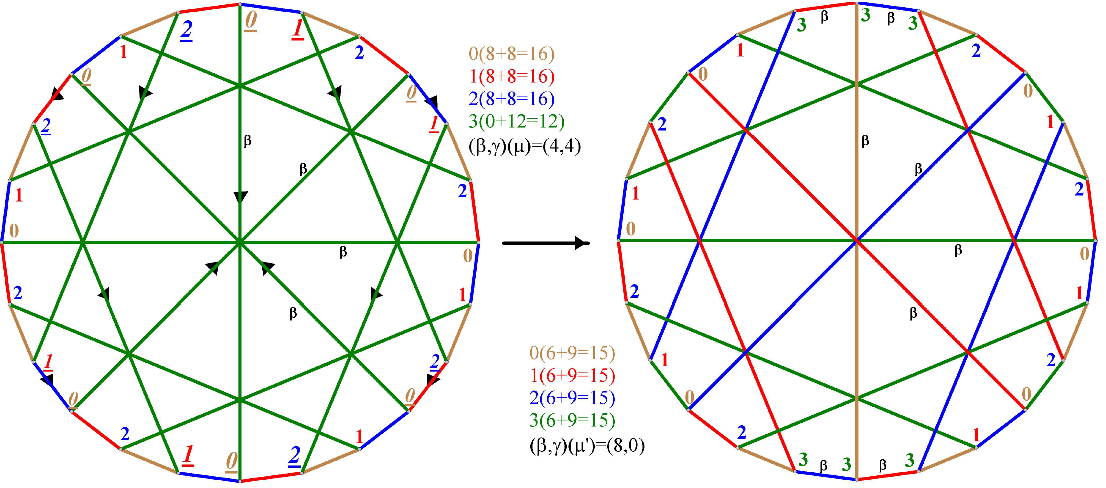}
\caption{From a lacunar sTC to an equitable sTC for the McGee cage (Example~\ref{mcgee}).}
\label{mcgill}
\end{figure}

\begin{example}\label{moka}
The M\"obius-Kantor graph $H(133)=H(10000101_2)$ is a symmetric bipartite cubic 16-vertex graph that is the generalized Petersen graph $GP(8,3)$. A lacunar sTC $\mu$ of $H(133)=H(10000101_2)$ is depicted in the upper-left of Figure~\ref{moka-joe} is a lacunar sTC $\mu$ with $(\beta,\gamma)(\mu)=(2,3)$ that is $\gamma$-reduced in the upper-right of Figure~\ref{moka-joe} to an sTC $\mu'$ with $(\beta,\gamma)(\mu)=(2,2)$. This in turn is $\gamma$-reduced to an equitable sTC with $(\beta,\gamma)(\mu)=(3,0)$ on the lower-left of Figure~\ref{moka-joe}, that can be $\beta$-reduced to an sTC with $(\beta,\gamma)(\mu)=(1,2)$ as in the lower-right of Figure~\ref{moka-joe}.   
\end{example}

\begin{example}\label{dy}
The Dyck graph $F_{032}A$ \cite{Foster} is the unique cubic symmetric graph on 32 vertices. $F_{032}A$ is represented in Figure~\ref{figdyck} via a sequence of five $\beta$-$\gamma$-reductions of sTCs from a lacunar sTC $\mu$ through sTCs $\mu',\mu'',\mu''',\mu''''$ and finally into an EqTC $\mu'''''$. 
\end{example}

\begin{figure}[htp]
\includegraphics[scale=0.55]{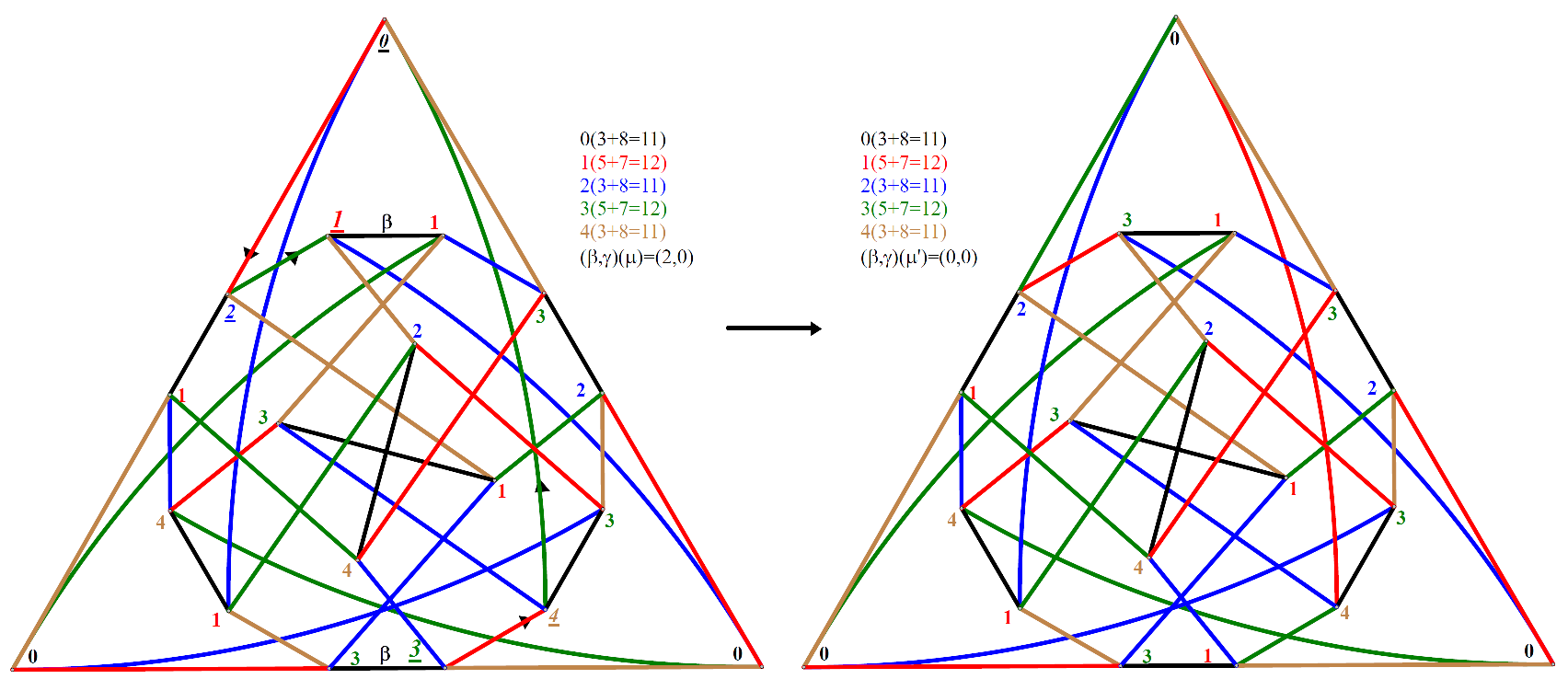}
\caption{Application to the Robertson 19-vertex cage (Example~\ref{petrob}).}
\label{figrob}
\end{figure}

\begin{example}\label{cansado} 
The distance transitive 90-vertex Foster graph $Fos$ with intersection array $\{3,2,2,2,2,1,1,1;1,1,1,1,2,2,2,3\}$ is hamiltonian with LCF notation $[17,-9,37,-37,9,$ $-17]^{15}$ seen in Figure~\ref{fig90} as the union of an external Hamilton 90-cycle $\eta$ and an internal 1-factor $F$ of 45 edges. 
 \end{example}
 
 \begin{theorem}\label{fost} There exists a 1-total-perfect lacunar sTC $\mu$ of the Foster graph $Fos$.\end{theorem}
 
 \begin{proof}
 $Fos$ is represented on the left of Figure~\ref{fig90} with a lacunar sTC $\mu$ that uses color set $\chi(3)$. 
There, $\eta$ is seen counterclockwise from the rightmost vertex, numbered 0 (and then 1, 2, etc., where not all numbers  $< 90$ are explicitly shown), via the color cycle $\mu(\eta)=((1_20_12_0)^{30})$, with edge colors appearing as subindices of their preceding vertex colors. In addition, $F$ shows only green (color 3) edges. The representation has 15 $\beta$-edges, so $\beta(Fos)\le15$. 
Also, the difference between the maximum and minimum cardinalities of the $\mu$-classes is 15, so $\gamma(\mu)=15$. 
 The hazel (0) $\mu$-class,
formed by the vertices bearing numbers congruent to 1 (mod 3), 
 induces a set of 15 green $\beta$ edges forming a total-perfect code. 
 \end{proof}

To obtain an equitable sTC on $Fos$ via a $\gamma$-reduction from $\mu$, consider the cycles mod 90 obtained by iterating the periodic additive operation pattern $(-1,+17,-1,+9;+1,+37,+1,\\ +9;\!)$, where each comma or semicolon is used as separator of the terms whose difference is the preceding signed-number. There are 3 such cycles, namely:

\begin{eqnarray}\label{tres}\begin{array}{l}
_{(1,0,17,16;\,25,26,63,64;\,73,72,89,88;\,7,8,45,46;\,55,54,71,70;\,79,80,27,28;\,37,36,53,52;\,61,62,9,10;\,19,18,35,34;\,43,44,81,82;)}\\
_{(31,30,47,46;\,55,56,3,4;\,13,12,29,28;\,37,38,75,76;\,85,84,11,10;\,19,20,57,58;\,67,66,83,82;\,1,2,39,40;\,49,48,65,64;\,73,74,21,22;)}\\
_{(61,60,77,76;\,85,86,33,34;\,43,42,59,58;\,67,68,15,16;\,25,24,41,40;\,49,50,87;\,7,6,23,22;\,31,32,69,70;\,79,78,5,4;\,13,14,51,52;)}\end{array}\end{eqnarray}

\noindent Each 6-subsequence formed by pairwise contiguous entries in one such 40-cycle constituted by three entries previous and three entries next to a semicolon represents an MCAP by just setting successively the corresponding vertex positions. Such 6-subsequences starting and ending in a number congruent to 0 or 2 mod 3 is assigned in Figure~\ref{fig90} color red (1) or blue (2), respectively. Each pair of contiguous terms in a cycle separated by a semicolon represents an MCAP formed by a green edge closing a 9-cycle with 8 successive edges of $\eta$.

\begin{figure}[htp]
\hspace*{2.5cm}
\includegraphics[scale=0.8]{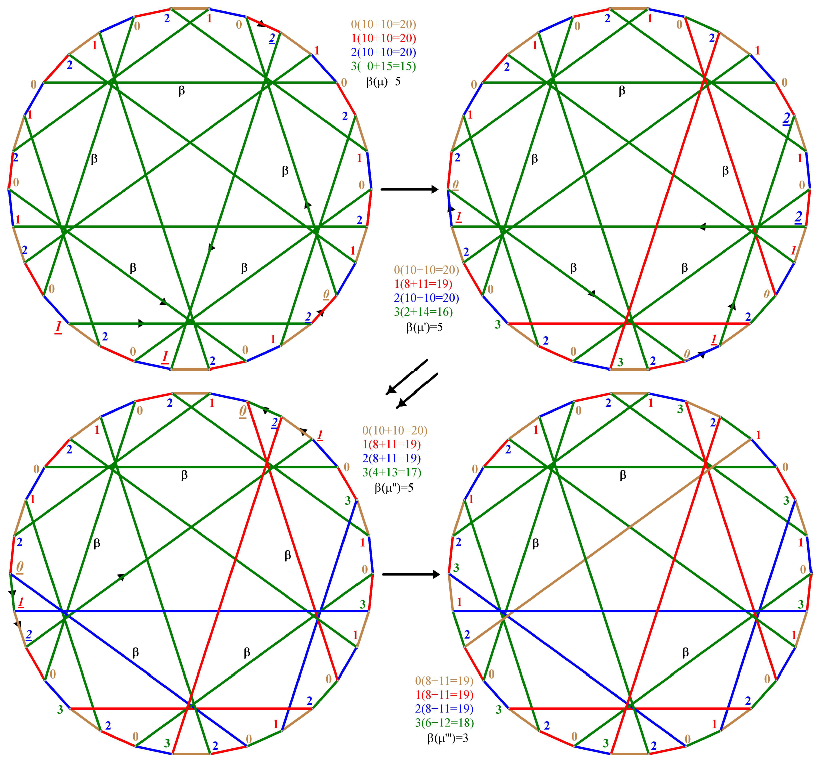}
\caption{From sTC $\mu$ of $TuCo$, $\beta(\mu)=5$, to equitable sTC $\mu'''$, $\beta(\mu''')=3$ (Example~\ref{higo}).}
\label{figtuco}
\end{figure}

An equitable sTC of $Fos$ is obtained as on the right of Figure~\ref{fig90}, where 12 $\gamma$-reductions on 12 corresponding MCAPs, where green (3) is exchanged with red (1), blue (2) and hazel (0) via four MCAPs each extracted from the first two 40-cycles in display (\ref{tres}), as shown on the expression of those 40-cycles on the bottom of Figure~\ref{fig90} (a pair of those 40-cycles expressed in two different ways), where red, blue and hazel fonts indicate (against the remaining black entries) the vertex numbers employed in those MCAPs. In the two representations of $Fos$, those vertex numbers $\neq 0$ and congruent to 0 mod 9 are shown in a larger type and all of them but 63 are in red font, indicating they are the eight endvertices of four MCAPs $(72,89,88;7,8,45)$, $(54,71,70;79,80,27)$, $(36,53,52;61,62,9)$ and $(18,35,34;43,44,81)$. On the other hand, the blue numbers indicate the eight endvertices of the four MCAPs whose numbers are congruent to 2 mod 9, (excluding 83), namely $(74,21,22;31,30,47)$, $(56,3,4;13,12,29)$, $(38,75,76;85,84,11)$ and $(2,39,40;49,48,65)$. Those sixteen endvertices and the eight endvertices of the four (length 1) remaining MCAPs, namely $(46;55)$, $(28;37)$, $(10;19)$ and $(64;73)$, are exchanged from their colors $\neq 3$ to green color, creating 14 new $\beta$-edges, and increasing the original $\mu(Fos)=15$ to a $\mu'(Fos)=29$, where $\mu'$ is the resulting equitable sTC. Observe that the selected set of 12 MCAPs 
insures distinct colors on adjacent edges and excludes a vertex color from its incident edges. Thus, Definition~\ref{acanomas} insures the sTC condition.  
Noticing that the listing of $\mu$ is $i(30+30=60)$, for $i=0,1,2$, and $3(0+45=45)$, the listing for $\mu'$ is $i(22+34=56)$, for $i=0,1,2$, and $3(24+33=57)$. 

\begin{example}\label{102} 
The distance transitive 102-vertex Biggs-Smith graph $BS$ \cite{BS} with intersection array $[3,2,2,2,1,1,1;1,1,1,1,1,1,3]$ is hamiltonian with LCF notation\\\,
[16, 24, -38, 17, 34, 48, -19, 41, -35, 47, -20, 34, -36, 21, 14, 48, -16, -36, -43, 28, -17, 21, 29, -43, 46, -24, 28, -38, -14, -50, -45, 21, 8, 27, -21, 20, -37, 39, -34, -44, -8, 38, -21, 25, 15, -34, 18, -28, -41, 36, 8, -29, -21, -48, -28, -20, -47, 14, -8, -15, -27, 38, 24, -48, -18, 25, 38, 31, -25, 24, -46, -14, 28, 11, 21, 35, -39, 43, 36, -38, 14, 50, 43, 36, -11, -36, -24, 45, 8, 19, -25, 38, 20, -24, -14, -21, -8, 44, -31, -38, -28, 37]\,
seen in Figure~\ref{fig102} as the union of an external Hamilton 102-cycle $\eta$ and an internal 1-factor $F$ of 51 edges. $BS$ is represented on the left of Figure~\ref{fig102} with a lacunar sTC $\mu$ that uses color set $\chi(3)$.  There, $\eta$ is seen counterclockwise from the rightmost vertex via the color cycle $\mu(\eta)=((1_02_10_2)^{34})$, where edge colors appear as subindices of their preceding vertex colors and $F$ has only green (color 3) edges. The representation has 17 $\beta$-edges, so $\beta(BS)\le 17$. 
Also, the difference between the maximum and minimum cardinalities of the $\mu$-classes is 17, so $\gamma(\mu)=17$. 
All vertices in $V(BG)=\{0,\ldots,,101\}$ divisible by 3 have red color, expressed in Figure~\ref{fig102} by citing their corresponding number. Then, the $\beta$-edges of $\mu$ are those green edges in thin trace in contrast to the remaining edges, in thick trace.

The LCF notation yields a list of the arcs of the 1-factor $F$ (two opposite arcs per edge) as in display (\ref{jodienda}), with arcs of $\beta$-edges expressed as $(x,y)$ and the remaining arcs as $[x,y]$. 
which is represented in Figure~\ref{fig102}. 

\begin{eqnarray}\label{jodienda}\begin{array}{llllllllll}
\,[0, 16]\!\!&\!\!(1, 25)\!\!&\!\![2, 66]\!\!&\!\![3, 20]\!\!&\!\![4,  38]\!\!&\!\!(5,  53)\!\!&\!\![6, 89]\!\!&\!\![7, 48]\!\!&\!\![8,  75]\!\!&\!\![9, 56]\!\!\\
\,[10,92]\!\!&\!\![11,45]\!\!&\!\!(12,78)\!\!&\!\!(13,34)\!\!&\!\![14,28]\!\!&\!\!(15,63)\!\!&\!\![16,0]\!\!&\!\!(17,83)\!\!&\!\![18,77]\!\!&\!\![19,47]\!\!\\
\,[20,3]\!\!&\!\!(21,42)\!\!&\!\![22,51]\!\!\!&\!\![23,82]\!\!&\!\![24,70]\!\!&\!\!(25,1)\!\!&\!\![26,54]\!\!&\!\![27,91]\!\!&\!\![28,14]\!\!&\!\![29,81]\\
(30,87)\!\!&\!\!(31,52)\!\!&\!\![32,40]\!\!&\!\!(33,60)\!\!&\!\!(34,13)\!\!&\!\![35,55]\!\!&\!\![36,101]\!\!&\!\!(37,76)\!\!&\!\![38,4]\!\!&\!\![39,97]\\
\,[40,32]\!\!&\!\![41,79]\!\!&\!\!(42,21)\!\!&\!\![43,68]\!\!&\!\!(44,59)\!\!&\!\![45,11]\!\!&\!\!(46,64)\!\!&\!\![47,19]\!\!&\!\![48,7]\!\!&\!\!(49,85)\\
\,[50,58]\!\!&\!\![51,22]\!\!&\!\!(52,31)\!\!&\!\!(53,5)\!\!&\!\![54,26]\!\!&\!\![55,35]\!\!&\!\![56,9]\!\!&\!\![57,71]\!\!&\!\![58,50]\!\!&\!\!(59,44)\\
(60,33)\!\!&\!\![61,99]\!\!&\!\!(62, 86)\!\!&\!\!(63,15)\!\!&\!\!(64,46)\!\!&\!\![65,90]\!\!&\!\![66,2]\!\!&\!\![67,98]\!\!&\!\![68,43]\!\!&\!\!(69,93)\\
\,[70,24]\!\!&\!\![71,57]\!\!&\!\![72,100]\!\!&\!\![73,84]\!\!&\!\!(74,95)\!\!&\!\![75,8]\!\!&\!\!(76,37)\!\!&\!\![77,18]\!\!&\!\!(78,12)\!\!&\!\![79,41]\\
\,[80,94]\!\!&\!\![81,29]\!\!&\!\![82,23]\!\!&\!\!(83,17)\!\!&\!\![84,73]\!\!&\!\!(85,49)\!\!&\!\!(86,62)\!\!&\!\!(87,30)\!\!&\!\![88,96]\!\!&\!\![89,6]\\
\,[90,65]\!\!&\!\![91,27]\!\!&\!\![92,10]\!\!&\!\!(93,69)\!\!&\!\![94,80]\!\!&\!\!(95,74)\!\!&\!\![96,88]\!\!&\!\![97,39]\!\!&\!\![98,67]\!\!&\!\![99,61]\\
\,[100,  72]\!\!&\!\![101,  36]\!\!&\!\!\!\!&\!\!\!\!&\!\!\!\!&\!\!\!\!&\!\!\!\!&\!\!\!\!&\!\!
\end{array}\end{eqnarray}

The $\beta$-edges have the following common endvertex colors, where a semicolon separates the vertex numbers of those 12 edges that change colors on the right of Figure~\ref{fig102}, in contrast to the remaining 5 commas:
$$\begin{array}{cl}
hazel:&((5;53),(17,83),(44;59),(62;86),(74;95);\\
red:&(12;78),(15;63),(21;42),(30;82),(33,60),(69,93);\\
blue:&(1;25),(13,34),(31;52),(37,76),(46,64),(49;85).\\
\end{array}$$

Those recolored 12 $\beta$-edges, (still in thin trace), 4 per each of the colors hazel, red and blue, have their endvertices recolored to green, corresponding to a total of 12 $\gamma$-reductions that takes the lacunar sTC $\mu$ with $\beta(\mu)=17$ on the left of Figure~\ref{fig102} to the equitable sTC $\mu'$ with $\beta(\mu')=23$ on its right.
 The listing of $\mu$ is  $I(34+34=68)$, for $i=0,1,2$, and $3(0+51=51)$. The listing for $\mu'$ is $i(26+38=64)$, for $i=0,1,2$, and $3(24+39=63)$.
 \end{example}

\begin{figure}[htp]
\hspace*{4mm}
\includegraphics[scale=1.25]{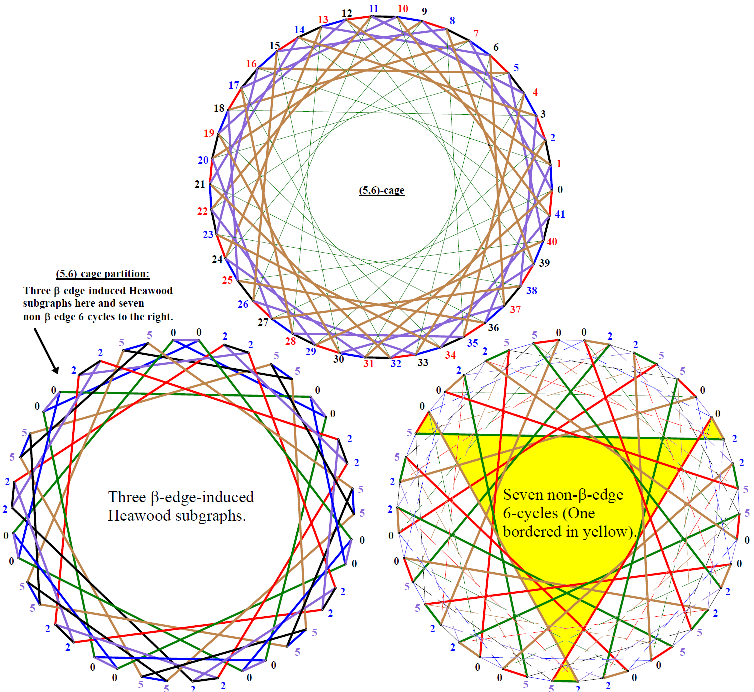}
\caption{Decomposition of the 42-vertex (5,6)-cage (Example~\ref{lastex}).}
\label{fig42ro}
\end{figure}

\section{Applications to cage graphs}\label{last}

As cases of cage graphs, we already treated $K_{3,3}$, the Heawood graph and the Petersen graph in Examples~\ref{K33}, ~\ref{e1} and \ref{Pap}, respectively, namely the $(3,4)$-,$(3,6)$- and $(3,5)$-cages on respective 6, 14 and 10 vertices. In addition, we consider the following cage graphs. 

\begin{example}\label{mcgee}
The McGee graph $McG$ is the only $(3,7)$-cage graph. $McG$ is hamiltonian with LCF notation $[-3,3]^4$, seen in Figure~\ref{mcgill} as the union of an external Hamilton 24-cycle $\eta$ and an internal 1-factor $F$ of 12 edges. $McG$ is represented on the left side of Figure~\ref{mcgill} with a lacunar sTC $\mu$ that uses color set $\chi(\Delta)=\chi(3)=\{0=$ hazel, $1=$ red, $2=$ blue, $3=$ green$\}$. This way, $\eta$ is seen clockwise from the top vertex via the color cycle
$\mu(\eta)=((0_21_02_1)^8)$, where edge colors appear as subindices of their preceding vertex colors, and $F$ has only green (color 3) edges.
There are four $\beta$-edges for $\mu$, distinguished via Greek letter $\beta$. By taking the four MCAPs constituted each solely by a $\beta$-edge and making corresponding $\gamma$-reductions, the equitable (non-lacunar) sTC $\mu'$ depicted on the right side of Figure~\ref{mcgill} is obtained, resumed with the notation $McG(15^4)$.
\end{example}

\begin{example}\label{petrob} We  apply Theorem~\ref{t1} to the Robertson graph $ Rob$ of 19 vertices, see Figure~\ref{figrob}.  $Rob$ is the $(4,5)$-cage.  
$Rob$ is shown on the left of Figure~\ref{figrob} with an sTC $\mu$ having just two $\beta$-edges. This is transformed into a TC $\mu'$ by exchanging the colors 3 (green) and 1 (red) along the MCAP 132103413 and making it into 312301431.   
\end{example}

\begin{example}\label{higo} The Tutte 8-cage \cite{GoRo} or Tutte-Coxeter graph $TuCo$  \cite{libro} is a cubic graph on 30 vertices and 45 edges. It is the unique $(3,8)$-cage graph and the unique $(3,8)$-Moore graph.
$TuCo$ is a vertex-transitive graph with intersection array $\{3,2,2,1;1,1,1,3\}$ and is 
hamiltonian with LCF notation $[-13,-9,7,-7,9,13]^5$ \cite{Frucht}. $TuCo$ is represented as the union of an external Hamilton 30-cycle $\eta$ and an internal 1-factor $F$ of 15 edges.
In the upper-left of Figure~\ref{figtuco} it is represented with a lacunar sTC $\mu$ that uses color set $\chi(\Delta)=\chi(3)$. 
This way, $\eta$ is seen clockwise from the lower-left vertex via the color cycle
$\mu(\eta)=((1_20_12_0)^{10})$, where edge colors appear as subindices of their preceding vertex colors, and $F$ has only green (color 3) edges.
There are five $\beta$-edges for $\mu$, distinguished via Greek letter $\beta$.

With terminal vertex $w_5$ being the one corresponding to the lower-left vertex color red (1) in $\mu$, an MCAP $S=w_0,e_1,w_1,e_2,w_2,e_3,w_3,e_4,w_4,e_5,w_5$ with color path $\mu(S)=1_32_10_30_12_31$ is distinguished, too, by having its vertex color numbers in underlined Italics in Figure~\ref{figtuco} (in contrast with the remaining vertex numbers, in Roman type), its two degree-one vertices in larger type and its edges oriented in the direction of $S$, so that $\mu(\rho(S))=\mu(w_0,e_1,e_2,e_3,e_4,e_5,w_5)=(1,3,1,3,1,3,1)$. 

\noindent The upper-right of Figure~\ref{figtuco} shows the resulting sTC $\mu'$  provided by Theorem~\ref{t1}. It was obtained by the ongoing $\gamma$-reduction. We note that $\mu(S)$ becomes on the right side of Figure~\ref{figtuco} into $\mu'(S)=3_12_3
0_10_32_13$. The listings for $\mu$ and $\mu'$ are shown near their representations, changing $i(10+10=20)$, for $i=0,1,2$, and $3(0+15=15)$ for $\mu$ into $i(10+10=20$, for $i=0,2$, $1(8+11=19)$ and  $3(2+14=16)$ for $\mu'$.
From $\mu'$ we $\gamma$-reduce to $\mu''$ via the MCAP $S'$ departing from the right-to-left middle horizontal arc in $F$ with color path $2_31_20_30_21_32$. This takes $\mu'$ into an sTC $\mu''$ in the lower-left of Figure~\ref{figtuco} with improved listing $0(10+10=20)$, $I(8+11=19)$, for $i=1,2$, and $3(4+13=17)$. From $\mu''$ we pass to an equitable sTC $\mu'''$ in the lower-right of Figure~\ref{figtuco} via an MCAP $S''$ starting at the leftmost vertex in $\eta$ with color path $0_31_02_31_02_30$. The listing for $\mu'''$ is
$i(8+11-19)$, for $i=0,1,2$, and $3(6+12=18)$. The last reduction works both as a total $\gamma$-reduction and as a partial $\beta$-reduction since $\beta(\mu''')=3<5=\beta(\mu'')=\beta(\mu')=\beta(\mu)$. 
\end{example}

\begin{example}\label{lastex} Consider the $(5,6)$-cage $\underline{\cap}_5^6$ on 42 vertices, depicted on the upper, lower-left and lower-right of Figure~\ref{fig42ro}. Then, $\underline{\cap}_5^6$ is hamiltonian with extended LCF notation $[(-11,-15,7),(11,15,-7)]^7$. In the top of Figure~\ref{fig42ro}, there is a lacunar semi-total coloring $\mu_{\underline{\cap}_5^6}$ of $\underline{\cap}_5^6$, with vertex periodic colors being just 0,1,2, selected from the palette of six colors available in this case, namely 0 = black, 1 = red, 2 = blue, 3 = green, 4 = hazel and 5 = violet. 
The vertices are numbered from 0 to 41 starting clockwise from the rightmost vertex, with colors 0,1,2 shown on the numbers by the corresponding visible colors black, red and blue. The edges of color 3 = green are the only $\beta$-edges, implying $\beta(\mu_{\underline{\cap}_5^6})\le 21$.

Another interesting semi-total coloring $\mu'_{\underline{\cap}_5^6}$ of ${\underline{\cap}_5^6}$ is shown with the (thick) edge set of ${\underline{\cap}_5^6}$ partitioned into those edges in the lower-left and those in the lower-right of Figure~\ref{fig42ro}, where the vertices are maintained in their position as in the upper depiction of ${\underline{\cap}_5^6}$. In the lower left, there is the disjoint union of three copies of the Heawood graph which are formed from all the $\beta$-edges of $\mu'_{\underline{\cap}_5^6}$. In the lower right, there is the disjoint union of seven copies of the 6-cycle $C_6$ with all the non-$\beta$-edges.  
\end{example}

\begin{theorem}\label{lasthm}
An sTC of the $(5,6)$-cage allows its decomposition into three copies of the Heawood graph via all the $\beta$-edges and seven 6-cycles via all the non-$\beta$-edges.
\end{theorem}

\begin{proof}
The claimed decomposition is obtained from the paragraph previous to the statement and Figure~\ref{fig42ro} as follows. The three claimed copies of the Heawood graph are formed by three corresponding Hamilton 14-cycles $H_0,H_1,H_2$ and respective complementary 1-factors $F_0,F_1,F_2$ presented in display (\ref{cage}) where vertex identification numbers are as established above and the second arrows in each row indicate the associated color 14-cycle or color 1-factor, where exponent 7 means color concatenation and multiplier 7 means the number of edges in each 1-factor.
\begin{eqnarray}\label{cage}\begin{array}{lllll}\hline
H_0&\rightarrow&(0,7,18,25,36,1,12,19,30,37,6,13,24,31)&\rightarrow&(5_05_4)^7\\
F_0&\rightarrow&(0,1)(6,7)(12,13)(18,19)(24,25)(30,31)(36,37)&\rightarrow& 7(5_25)\\\hline
H_1&\rightarrow&(2,9,20,27,38,3,14,21,32,39,8,15,26,33)&\rightarrow&(2_52_1)^7\\
F_1&\rightarrow&(2,3)(8,9)(14,15)(20,21)(26,27)(32,33)(38,39)&\rightarrow&7(2_02)\\\hline
H_2&\rightarrow&(4,11,22,29,40,5,16,23,34,41,10,17,28,35)&\rightarrow&(0_20_1)^7\\
F_2&\rightarrow&(4,\!5)(10,\!11)(16,\!17)(22,\!23)(27,\!28)(34,\!35)(40,\!41)&\rightarrow&7(0_50)\\\hline
\end{array}\end{eqnarray}
The seven claimed copies of the 6-cycle, say $C_i$ ($i=0,\ldots,6$), are presented in a similar fashion in display (\ref{siete}). 
\begin{eqnarray}\label{siete}\begin{array}{lllll}
C_0&\rightarrow&(1,16,15,30,29,2)&\rightarrow&(5_10_42_3)^2\\
C_1&\rightarrow&(7,22,21,36,35,8)&\rightarrow&(5_10_42_3)^2\\
C_2&\rightarrow&(13,28,27,0,41,14)&\rightarrow&(5_10_42_3)^2\\
C_3&\rightarrow&(19,34,33,6,5,20)&\rightarrow&(5_10_42_3)^2\\
C_4&\rightarrow&(25,40,39,12,11,26)&\rightarrow&(5_10_42_3)^2\\
C_5&\rightarrow&(31,4,3,18,17,32)&\rightarrow&(5_10_42_3)^2\;\mbox{ (bordered in yellow)}\\
C_6&\rightarrow&(37,10,9,24,23,38)&\rightarrow&(5_10_42_3)^2\\
\end{array}\end{eqnarray}
Both displays, (\ref{cage}) and (\ref{siete}), allow to establish the claimed decompositiont.
\end{proof}

\begin{question}
Are there any other semi-total colorings of regular graphs having the $\beta$-edges forming copies of interesting subgraphs, as in Example~\ref{lastex}?
\end{question}

\section{Applications to vertex expansions}\label{s3}

\begin{figure}[htp]
\includegraphics[scale=0.85]{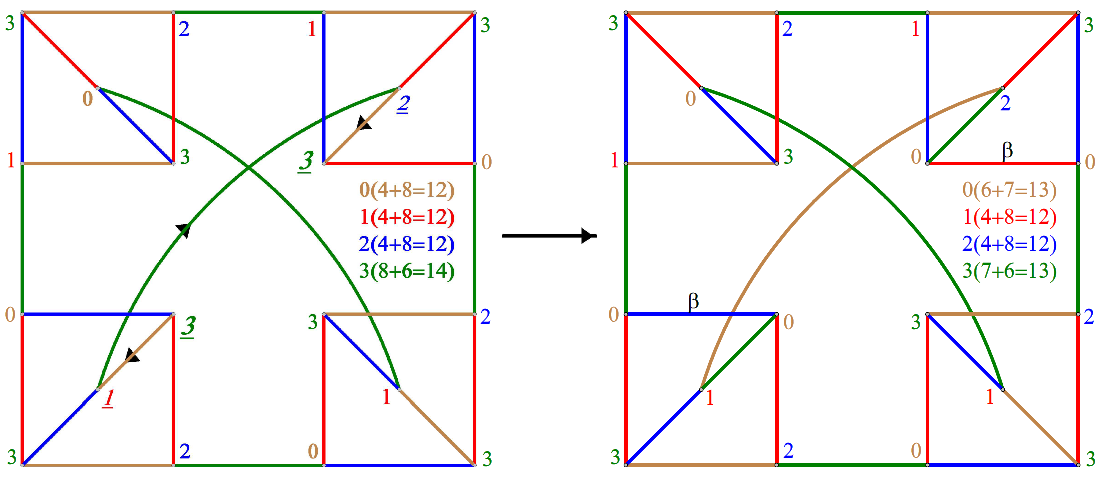}
\caption{From a TC of the graph $K_4^{K_{2,3}}$ to an equitable sTC (Example~\ref{brasil}).}
\label{sas1}
\end{figure}

\begin{figure}[htp]
\includegraphics[scale=0.88]{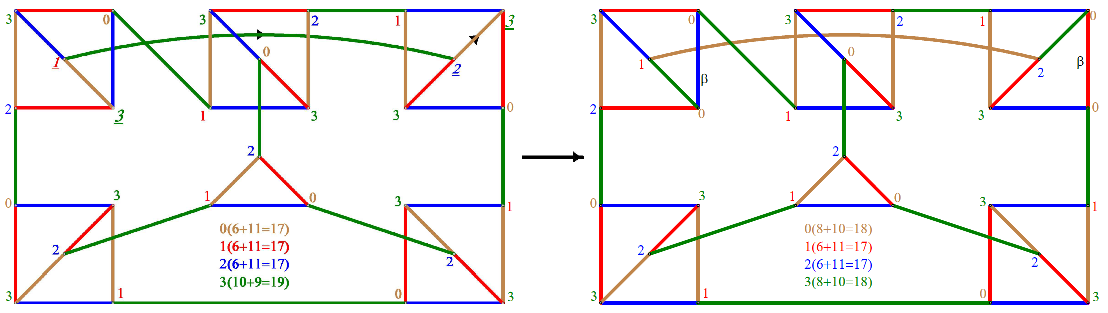}
\caption{From a TC of the graph $(K_3\square K_2)^{K_{2,3},K_3}$ to an equitable sTC (Example~\ref{brasil}).}
\label{sas2}
\end{figure}

\begin{example}\label{brasil}
Dantas et al. \cite{Dantas} show that for every cubic graph $G$, the cubic graph $G^{K_{2,3}}$ obtained from $G$ by replacing every vertex by a copy of $K_{2,3}$ is of type 1 and has equitable total chromatic number 5. They also show that by replacing one copy of $K_{2,3}$ in such a graph by a triangle $K_3$ the resulting graph $G^{K_{2,3},K_3}$ also has equitable total chromatic number 5 if $|V(G)|\ge 6$. 
 Two examples of TCs in such graphs of \cite{Dantas} are given in Figures~\ref{sas1} and~\ref{sas2}, showing how to $\gamma$-reduce them into equitable sTC's.

A TC $\mu:K_{2,3}\rightarrow[3]_0$ of the complete bipartite graph $K_{2,3}$ with vertex parts $\{u,v\}$ and $\{w,x,y\}$ will be denoted $\mu(K_{2,3}(u,v;w,x,y)=(a_{b'c'd'},a_{b''c'',d''};b,c,d)$, where $\{a,b,c,d\}=[3]_0$ and $a=\mu(u)=\mu(v)=\notin\{b,c,d\}=\{b',c',d'\}=\{b'',c'',d''\}$, with $\mu(w)=b,\mu(x)=c,\mu(y)=d$ and
$\mu(u,w).=b',\mu(u,x)=c',\mu(u,y)=d',\mu(v,w)-b'',\mu(v,x)=c'',\mu(v,y)=d''$.  

Figure~\ref{sas1} shows how to $\gamma$-reduce a TC $\mu$ of the graph $K_4^{K_{2.3}}$ to an equitable sTC. Numbering its vertices first left-to-right and then top to bottom from 0 to 19, it is seen on the left of the figure that $G$ is composed by copies of $K_{2,3}$ with
$\mu(K_{2,3}(0,7;1,4,6))=(3_{012},3_{120};2,0,1)$, $\mu(K_{2,3}(2,9;3,5,8))=(3_{012},3_{201};1,2,0)$, $\mu(K_{2,3}(10,17;11,14,16))=(3_{201},3_{120};0,1,2)$,\\ $\mu(K_{2,3}(13,18;12,15,19))=(3_{021},3_{102};2,1,0)$, and additional green (color 3) edges $(1,2)$, $(4,15)$, $(5,14)$, $(6,10)$, $(9,13)$, $(17,18)$.  
The MCAP $S=(11,(11,14),14,(14,6),6,(6,9),9)$ has skeleton $\rho(S)$ with color image
$3_01_32_03$, where edge colors appear as subindices of their preceding vertex colors, that is $\gamma$-reduced to $0_31_02_30$, as indicated by the horizontal arrow from the left color representation of the graph, with listing $i(4+8=12)$, for $i=0,1,2$, and $3(8+6=14)$, yielding $(\beta,\gamma)(\mu)=(0,2)$, to the right color representation, with listing $i(6+7=13)$, for $i=0,3$, and $j(4+8)$, for $j=1,2$, yielding $(\beta,\gamma)(\mu')=(2,1)$ with $\beta$-edges $(9,10)$ and $(11,12)$. 

In the following paragraph, a TC $\mu:K_3\rightarrow[3]_0$ with $V(K_3)=\{w,x,y\}$ will be denoted $\mu(K_3(w,x,y))=(a_cb_ac_b)$, where $\{a,b,c\}\subset[3]_0$ with $\mu(w)=a$, $\mu(x)=b$, $\mu(y)=c$ and
$\mu(w,x).=c$, $\mu(x,y)=a$ and $\mu(y,w)=b$.  

Figure~\ref{sas2} shows how to $\gamma$-reduce a TC $\mu$ of the graph $G=(K_3\square K_2)^{K_{2.3},K_3}$ to an equitable sTC. Numbering the vertices of $G$ first left-to-right and then top to bottom from 0 to 27, it is seen on the left of the figure that $G$ is composed by copies of $K_{2,3}$ with colorings
$\mu(K_{2,3}(0,10;1,6,9))=(3_{120},3_{201};0,1,2)$, $\mu(K_{2,3}(2,12;3,7,11))=(3_{120},3_{012};2,0,1)$, $\mu(K_{2,3}(5,13;4,8,14))=(3_{201},3_{012};1,2,0)$, $\mu(K_{2,3}(17,24;16,22,25))=(3_{210},3_{102};0,2,1)$,\\ $\mu(K_{2,3}(20,27;21,23,26))=(3_{201},3_{012};1,2,0)$, a copy of $K_3$ with coloring $\mu(K_3(15,18,19))=(2_01_20_1)$ and additional green (color 3) edges $(1,11)$, $(3,4)$, $(6,8)$, $(7,15)$, $(9,16)$, $(14,21)$, $(18,22)$, $(19,23)$, $(25,26)$.
The MCAP $S=(6,(6,8),8,(8,5),5)$ has skeleton $\rho(S)$ with color image
$1_32_03$, where edge colors appear as subindices of their preceding vertex colors, that is $\gamma$-reduced to $1_02_30$, as indicated by the horizontal arrow from the left color representation of the graph, with listing $i(6+11=17)$, for $i=0,1,2$, and $3(10+9=19)$, yielding $(\beta,\gamma)(\mu)=(0,2)$, to the right color representation, with listing $i(8+10=18)$, for $i=0,3$, and $j(6+11)$, for $j=1,2$, yielding $(\beta,\gamma)(\mu')=(2,1)$ with $\beta$-edges $(1,10)$ and $(5,12)$. 
\end{example}

\end{document}